\pdfoutput=1
\RequirePackage{ifpdf}
\ifpdf 
\documentclass[pdftex]{sigma}
\else
\documentclass{sigma}
\fi

\usepackage{mathrsfs}
\numberwithin{equation}{section}
\newtheorem{Theorem}{Theorem}[section]
\newtheorem{Corollary}[Theorem]{Corollary}
\newtheorem{Lemma}[Theorem]{Lemma}
\newtheorem{Proposition}[Theorem]{Proposition}
{\theoremstyle{definition}
\newtheorem{Definition}[Theorem]{Definition}
\newtheorem{Remark}[Theorem]{Remark}
}

\DeclareMathOperator{\Aut}{Aut}
\DeclareMathOperator{\Image}{Image}
\DeclareMathOperator{\Sing}{Sing}
\DeclareMathOperator{\id}{id}

\begin{document}

\allowdisplaybreaks

\renewcommand{\thefootnote}{$\star$}

\renewcommand{\PaperNumber}{046}

\FirstPageHeading

\ShortArticleName{Scalar Flat K\"ahler Metrics on Af\/f\/ine Bundles over $\mathbb{CP}^1$}

\ArticleName{Scalar Flat K\"ahler Metrics\\ on Af\/f\/ine Bundles over $\boldsymbol{\mathbb{CP}^1}$\footnote{This paper is
a~contribution to the Special Issue on Progress in Twistor Theory.
The full collection is available at
\href{http://www.emis.de/journals/SIGMA/twistors.html}{http://www.emis.de/journals/SIGMA/twistors.html}}}

\Author{Nobuhiro HONDA}

\AuthorNameForHeading{N.~Honda}

\Address{Mathematical Institute, Tohoku University, Sendai, Miyagi, Japan}
\Email{\href{mailto:honda@math.tohoku.ac.jp}{honda@math.tohoku.ac.jp}}

\ArticleDates{Received November 12, 2013, in f\/inal form April 15, 2014; Published online April 19, 2014}

\Abstract{We show that the total space of any af\/f\/ine $\mathbb{C}$-bundle over $\mathbb{CP}^1$ with negative degree
admits an ALE scalar-f\/lat K\"ahler metric.
Here the degree of an af\/f\/ine bundle means the negative of the self-intersection number of the section at inf\/inity
in a~natural compactif\/ication of the bundle, and so for line bundles it agrees with the usual notion of the degree.}


\Classification{53A30}

\renewcommand{\thefootnote}{\arabic{footnote}}
\setcounter{footnote}{0}

\section{Introduction}

If $\Gamma$ is a~f\/inite subgroup of the unitary group ${\rm{U}}(2)$ which acts freely on the unit sphere around the
origin in $\mathbb{C}^2$, it is natural to ask existence of a~K\"ahler metric def\/ined on the minimal resolution of the
quotient space $\mathbb{C}^2/\Gamma$, which has `small curvature' and which is asymptotically locally Euclidean (ALE) at
inf\/inity.
If $\Gamma$ is a~f\/inite subgroup of ${\rm{SU}}(2)$, Kronheimer~\cite{Kr89a} constructed ALE Ricci-f\/lat K\"ahler
metrics on the minimal resolution by means of so called the hyperK\"ahler quotient, and further showed~\cite{Kr89b} that
the metrics are determined by the the cohomology classes of a~collection of K\"ahler forms associated to the
hyper-K\"ahler structure.
When $\Gamma$ is a~f\/inite cyclic subgroup of ${\rm{U}}(2)$ generated by scalar matrices, $\Gamma$ is not included in
${\rm{SU}}(2)$ unless $|\Gamma|=2$, and the minimal resolution of $\mathbb{C}^2/\Gamma$ is simply the total space of the
line bundle $\mathscr O(-n)$, where $n=|\Gamma|$, and the unique negative section is the exceptional locus of the
resolution.
Because the section intersects positively with the canonical class if $n>2$, there exists no Ricci-f\/lat K\"ahler
metric on $\mathscr O(-n)$ if $n>2$.
LeBrun~\cite{LB88} constructed on this complex surface a~scalar-f\/lat K\"ahler (SFK) metric which is also ALE.
The metric is invariant under a~natural ${\rm{U}}(2)$-action, and may be considered to be the natural K\"ahler metric on
$\mathscr O(-n)$.
Later, Calderbank--Singer~\cite{CS2004} pointed out that, for any cyclic subgroup $\Gamma\subset {\rm{U}}(2)$, the
minimal resolution of $\mathbb{C}^2/\Gamma$ admits an ALE SFK metric.
All these metrics are anti-self-dual (ASD) with respect to the complex orientation.

ALE spaces can be compactif\/ied to be an orbifold by adding a~point at inf\/inity, and after an appropriate conformal
change the metric can be extended to the compactif\/ication as an ASD metric on the orbifold.
Small deformations of ASD conformal structures are governed by a~deformation complex, and if the space is compact, the
index of the complex is expressed in terms of topological invariants of the space.
Viaclovsky~\cite{V} computed the index of the deformation complex for various compact orbifolds in explicit form, and
show in particular that if the space is the compactif\/ication $\widehat{\mathscr O(-n)}$ of $\mathscr O(-n)$, the index
is $12-4n$.
As the obstruction for the deformation complex vanishes, this means that LeBrun's metric has a~non-trivial deformation
as an ALE ASD metric (if $n>3$).
In~\cite{HonCMP2} we computed the ${\rm{U}}(2)$-action on the relevant cohomology group in concrete form, and computed
the dimension of the moduli space of ALE ASD metrics near the LeBrun metric.
Also we found that there exists a~real 1-parameter family of deformation of the LeBrun metric which preserves not only
ASD-ALE property but also a~K\"ahler representative.

In this article we investigate all small deformations of the LeBrun metric on $\mathscr O(-n)$ {\textit{as an ALE SFK metric}}.
If the complex structure is f\/ixed, the following rigidity is shown:
\begin{Proposition}[=~Proposition~\ref{prop:rigid}]
\label{prop:rigid0}
When we fix the complex structure on $\mathscr O(-n)$, LeBrun's metric on $\mathscr O(-n) $ cannot be deformed as an
ALE SFK metric by small deformations.
\end{Proposition}

In order to explain what happens when we allow the complex structure on $\mathscr O(-n)$ to vary in deformations, we
recall that the line bundle $\mathscr O(-n)$ is included as a~special member of an $(n-1)$-parameters family of
af\/f\/ine $\mathbb{C}$-bundles over $\mathbb{CP}^1$, whose transition law for f\/iber coordinates $\zeta_0$ over
$U_0=\mathbb{C}(u)\subset\mathbb{CP}^1$ and $\zeta_1$ over $U_1=\mathbb{C}(1/u)\subset\mathbb{CP}^1$ is concretely given
by
\begin{gather}
\label{af00}
\zeta_0 = \frac 1{u^n} \zeta_1 + \sum\limits_{l=1}^{n-1} \frac {t_l} {u^l},
\qquad
(t_1,\dots,t_{n-1})\in \mathbb{C}^{n-1},
\end{gather}
where $t_l$-s are parameters.
If $(t_1,\dots,t_{n-1})=(0,\dots,0)$, this gives the line bundle $\mathscr O(-n)$,
but if $(t_1,\dots,t_{n-1})\neq(0,\dots,0)$, the af\/f\/ine bundle~\eqref{af00} has no global section and it is just an af\/f\/ine bundle.
For $t=(t_1,\dots,t_{n-1})\in\mathbb{C}^{n-1}$ we denote by $A_t$ for the af\/f\/ine $\mathbb{C}$-bundle over
$\mathbb{CP}^1$ def\/ined by~\eqref{af00}.
Then we prove the following

\begin{Theorem}
\label{thm:0}
There exists a~neighborhood $B\subset\mathbb{C}^{n-1}$ of the origin for which LeBrun's ALE SFK metric on $\mathscr
O(-n)$ extends naturally to $A_t$ if $t\in B$, as ALE SFK metrics.
\end{Theorem}

This will be shown as Theorem~\ref{thm:main01}, and from the proof, this family of metrics can be regarded as the versal
family for the LeBrun metric on $\mathscr O(-n)$ as ALE SFK metrics.
The 1-parameter family of ALE SFK metrics on the 4-manifold $\mathscr O(-n)$ obtained in~\cite{HonCMP2} is exactly the
restriction of the family of ALE SFK metrics in Theorem~\ref{thm:0} to the f\/irst (or the last) coordinate axis.

Next for explaining an immediate consequence of Theorem~\ref{thm:0}, we recall that any af\/f\/ine $\mathbb{C}$-bundle
over $\mathbb{CP}^1$ can be naturally compactif\/ied to a~Hirzebruch surface by attaching a~section at inf\/inity.
We call the negative of the last self-intersection number as \textit{degree} of the af\/f\/ine bundle.
Then for any $t\in\mathbb{C}^{n-1}$ the af\/f\/ine bundle def\/ined by the transition law~\eqref{af00} is of degree~$n$.
Conversely, if $n>1$, any af\/f\/ine $\mathbb{C}$-bundle over $\mathbb{CP}^1$ of degree $n$ is of the form $A_t$ for
some $t\in\mathbb{C}^{n-1}$.
Now because the equation~\eqref{af00} is linear in the variables $\zeta_0$, $\zeta_1$, $t_1, t_2,\dots,t_{n-1}$, the af\/f\/ine
bundle $A_t$ and $A_{ct}$ is isomorphic for any $c\in \mathbb{C}^*$.
Therefore, Theorem~\ref{thm:0} implies the following

\begin{Corollary}[=~Corollary~\ref{cor:exaf}]
Any affine $\mathbb{C}$-bundle over $\mathbb{CP}^1$ of negative degree $($in the above sense$)$ admits an ALE SFK metric.
\end{Corollary}

Finally we explain some property of the family of ALE SFK metrics obtained in Theorem~\ref{thm:0}.
In contrast with the LeBrun metric for which the rigidity holds as in Proposition~\ref{prop:rigid0}, for the deformed
metrics, we have the following

\begin{Proposition}\label{prop:var}
Let $B\subset\mathbb{C}^{n-1}$ be as in Theorem~{\rm \ref{thm:0}}, and for $t\in B$ let $g_t$ be the ALE SFK metric on~$A_t$.
Then if $t\neq 0$ and $t$ is sufficiently close to the origin, there exists a~smooth arc $\gamma_t \subset B$
passing through the point $t$ which satisfies the following:
\begin{itemize}\itemsep=0pt
\item[$(i)$] the complex structure of $A_t$ is constant along the arc $\gamma_t$,
\item[$(ii)$] the conformal class of the ALE SFK metric $g_t$ varies when $t$ moves along $\gamma_t$.
\end{itemize}
\end{Proposition}

\section{Preliminary computations for Hirzebruch surfaces}
\label{s:Hirz}

\subsection{Notation and convention}
\label{ss:notation}

For an integer $n\ge 0$, $\mathbb{F}_n$ denotes the Hirzebruch surface of degree $n$; namely $\mathbb{F}_n=\mathbb
P(\mathscr O(-n)\oplus \mathscr O)$ over $\mathbb{CP}^1$.
We write $\pi:\mathbb{F}_n\to\mathbb{CP}^1$ for the projection, and $f$ for a~f\/iber (class) of $\pi$.
We denote $\Gamma_0$ for $(-n)$-section of $\pi$, which is unique when $n>0$.
We have $H^2(\mathbb{F}_n,\mathbb{Z})\simeq {\rm{Pic}}\mathbb{F}_n\simeq \mathbb{Z}[\Gamma_0]\oplus\mathbb{Z}[f]$, and
$-K_{\mathbb{F}_n}\simeq \mathscr O(2\Gamma_0 + (n+2)f)$ for the anticanonical class.
${\rm{Aut}}_0\mathbb{F}_n$ denotes the identity component of holomorphic transformation group of $\mathbb{F}_n$, and
for a~section~$L$ of $\mathbb{F}_n\to\mathbb{CP}^1$, ${\rm{Aut}}_0(\mathbb{F}_n,L)$ denotes the subgroup of
${\rm{Aut}}_0\mathbb{F}_n$ consisting of transformations which keep~$L$ invariant.
If $n>0$ and $L=\Gamma_0$, we have ${\rm{Aut}}_0\mathbb{F}_n= {\rm{Aut}}_0(\mathbb{F}_n,L)$.
Two pairs $(\mathbb{F}_n,L)$ and $(\mathbb{F}_n,L')$ are called isomorphic as a~pair if there is a~biholomorphic map
$\phi:\mathbb{F}_n \to \mathbb{F}_n$ which satisf\/ies $\phi(L) = L'$.
$\Gamma_{\infty}$ means a~section whose self-intersection number is $(+n)$.
${\rm{Aut}}_0\mathbb{F}_n$ acts transitively on the space of $(+n)$-sections, and $\Gamma_{\infty}$ may be
identif\/ied with the section $\mathbb P(\mathscr O(-n))$.
Thus the complement $\mathbb{F}_n\backslash \Gamma_{\infty}$ can be identif\/ied with the total space of the line bundle~$\mathscr O(-n)$.

We write the linear system to which a~section of $\mathbb{F}_n\to\mathbb{CP}^1$ belongs in the form $|\Gamma_0+kf|$, $k\ge 0$.
It is well-known that this system has an irreducible member only when $k=0$ or $k\ge n$.
So if~$L$ is a~section with positive self-intersection number, we have~$L\in |\Gamma_0 + (n+l)f|$ for some $l\ge 0$.
The letter $l$ is always used in this meaning throughout the article.
We have $\Gamma_{\infty}\in |\Gamma_0+nf|$.
The system $|\Gamma_0 + (n+l)f|$ is very ample if and only if $l>0$.
Moreover we have $h^0(\mathscr O_{\mathbb{F}_n}(\Gamma_0+(n+l)f))= n+2l+2$, where $h^0$ means $\dim H^0$.
Therefore the complement of any member of these systems is realized in an af\/f\/ine space $\mathbb{C}^{n+2l+1}$.

We will also use the following result regarding the dimension of the cohomology group $H^i(\Theta_{\mathbb{F}_n})$ of
the tangent sheaf.
Namely if $n>0$ we have \cite{MKbook}
\begin{gather}
\label{cohohir}
h^0(\Theta_{\mathbb{F}_n}) = n + 5,
\qquad
h^1(\Theta_{\mathbb{F}_n}) = n - 1,
\qquad
h^2(\Theta_{\mathbb{F}_n}) = 0.
\end{gather}
This will also be shown in the proof of Proposition~\ref{prop:l0}.

\subsection[Af\/f\/ine $\mathbb{C}$-bundles over $\mathbb{CP}^1$]{Af\/f\/ine
$\boldsymbol{\mathbb{C}}$-bundles over $\boldsymbol{\mathbb{CP}^1}$}
\label{ss:ab}

Let ${\rm Af}(\mathbb{C})$ be the group of complex af\/f\/ine transformations of $\mathbb{C}$; namely those of the form
$\zeta\mapsto a\zeta + b$ for $\zeta\in \mathbb{C}$, where $a\in \mathbb{C}^*$ and $b\in \mathbb{C}$.
Let $X$ be a~projective algebraic manifold.
By an af\/f\/ine $\mathbb{C}$-bundle over~$X$, we mean (as usual) a~$\mathbb{C}$-bundle $A\to X$ whose structure group
is ${\rm {Af}}(\mathbb{C})$.
In this subsection, according to~\cite{At55}, we brief\/ly explain a~classif\/ication of af\/f\/ine $\mathbb{C}$-bundles
over~$X$, and then apply it to a~concrete description of af\/f\/ine $\mathbb{C}$-bundles over~$\mathbb{CP}^1$.

As in the case of any f\/iber bundle with prescribed structure group, isomorphic classes of af\/f\/ine
$\mathbb{C}$-bundles over $X$ are naturally in 1-1 correspondence with the cohomology set $H^1(X, \mathscr {A}f)$, where
$\mathscr {A}f$ means the sheaf of germs of holomorphic maps from open sets in~$X$ to the group~${\rm{Af}}(\mathbb{C})$.
The set $H^1(X, \mathscr {A}f)$ is of course the inductive limit of $H^1(\mathscr U,\mathscr Af)$ with respect to open
covering $\mathscr U$-s of~$X$.
For each $\mathscr U$, there is a~natural map
\begin{gather*}
\rho_{\mathscr U}: \ H^1(\mathscr U,{\mathscr A}f)
 \longrightarrow
H^1(\mathscr U,\mathscr O^*),
\end{gather*}
which is induced from the natural homomorphism $ {\rm {Af}}(\mathbb{C})\to \mathbb{C}^*$ that takes the coef\/f\/icient
of the linear part.
These naturally induce a~map $\rho:H^1(X, \mathscr Af)\to H^1(X,\mathscr O^*)$.
Therefore we have
\begin{gather}
\label{star1}
H^1(X,\mathscr Af) \simeq \bigsqcup_{\xi \in H^1(X,\mathscr O^*)}\rho^{-1}(\xi).
\end{gather}
Geometrically, for an (isomorphism class of) af\/f\/ine bundle $A\in H^1(X,\mathscr Af)$, the image $\rho(A)\in
H^1(X,\mathscr O^*)$ is exactly (the isomorphism class of) the dual line bundle of the normal bundle of the section
$\overline A \backslash A$ in $\overline A$, where $\overline A$ means the compactif\/ied $\mathbb{CP}^1$-bundle which
is obtained from $A\to X$ by the standard inclusion ${\rm Af}(\mathbb{C})\subset\mathrm{PGL}(2,\mathbb{C})$.

Returning to the \v Cech cohomology group, analogously to~\eqref{star1}, we clearly have, for each open covering
$\mathscr U$ of $X$,
\begin{gather*}
H^1 (\mathscr U,\mathscr Af) = \bigsqcup_{\xi \in H^1(\mathscr U,\mathscr O^*)}\rho_{\mathscr U}^{-1}(\xi).
\end{gather*}
In order to describe the set $\rho_{\mathscr U}^{-1}(\xi)$, we write $\mathscr U = \{U_i\}$, and let $\xi\in
H^1(\mathscr U,\mathscr O^*)$ be represented by a~1-cocycle $\{a_{ij}\}$, so that $a_{ij}\in \mathscr O^*(U_{ij})$ where
$U_{ij}=U_i\cap U_j$.
We f\/ix a~collection $\{h_i\}$ of meromorphic functions, where $h_i$ is def\/ined on $U_i$, that satisfy $h_i = a_{ij}
h_j$ on $U_{ij}$ (this is possible from the projectivity assumption for $X$), and let $D$ be the divisor def\/ined by
$\{h_i = 0\}$.
Though this is not necessarily ef\/fective, it is `linear equivalent' to the line bundle $\xi$.
Under these f\/ixing of $\{a_{ij}\}$ and $\{h_i\}$, let $\{(a'_{ij},b'_{ij})\}\in H^1(\mathscr U,\mathscr Af)$, where
$a'_{ij}\in\mathscr O^*(U_{ij})$ and $b'_{ij}\in \mathscr O(U_{ij})$, be a~representative of an element of
$\rho_{\mathscr U}^{-1}(\xi)$.
If we choose any $\phi_i\in \mathscr O^*(U_i)$ and $\psi_i\in \mathscr O(U_i)$ for each~$i$ and apply a~f\/iber
coordinate change $\tilde{\zeta}_i= \phi_i\zeta_i + \psi_i$ on $U_i$, then the new 1-cocycle $\{(\tilde a_{ij},\tilde
b_{ij})\}$ associated to $\{\tilde{\zeta}_i\}$, which is another representative of the same element of $H^1(\mathscr
U,\mathscr Af)$, is readily seen to be given by
\begin{gather}
\label{cd88}
\tilde a_{ij} = \frac{\phi_i}{\phi_j}a'_{ij},
\qquad
\tilde b_{ij} = \phi_i b'_{ij} + \psi_i - \frac{\phi_i}{\phi_j}\psi_j a'_{ij}.
\end{gather}
The f\/irst equation of these shows that any element of $\rho_{\mathscr U}^{-1} (\xi)$ can be represented by a~cocycle
of the form $\{(a_{ij},b'_{ij})\}$ (namely, by using the original representative $\{a_{ij}\}$ for $\xi\in H^1(\mathscr
U,\mathscr O^*)$), and in the following, for any element of $\rho_{\mathscr U}^{-1} (\xi)$, we only consider such
representatives.
This means that we only consider f\/iber coordinate changes $\{(\phi_i,\psi_i)\}$ which satisfy $\phi_i = \phi_j$ on
$U_{ij}$, and hence we can write $\phi_i = t$ for all~$i$ for some constant $t\in \mathbb{C}^*$.
Then the second equation of~\eqref{cd88} becomes (after replacing $a'_{ij}$ by $a_{ij}$)
\begin{gather}
\label{cc89}
\tilde b_{ij} = t b'_{ij} + \psi_i -\psi_j a_{ij}.
\end{gather}
This is the transformation law for representatives, under coordinate changes that satisfy the above constraint.

We are still f\/ixing $\mathscr U=\{U_i\}$, $\xi \in H^1(\mathscr U,\mathscr O^*)$, a~representative $\{a_{ij}\}$ of~$\xi$,
and $\{h_i\}$ that satisf\/ies $h_i=a_{ij}h_j$ on $U_{ij}$.
If $\{(a_{ij},b_{ij})\}$ is a~1-cocycle that represents an element of $\rho^{-1}_{\mathscr U}(\xi)$, we def\/ine
\begin{gather*}
c_{ij}:= \frac{b_{ij}}{h_i}
\qquad
\text{on}
\quad
U_{ij}.
\end{gather*}
Then from the cocycle condition for $\{(a_{ij},b_{ij})\}$, it follows that $\{c_{ij}\}$ is a~1-cocycle whose value is in
$\mathscr O(D)$, where $\mathscr O(D)$ is the sheaf of holomorphic functions $f$ for which $fh_i$ is holomorphic for any~$i$.
Moreover if we apply f\/iber coordinate changes of the form $\{(\phi_i,\psi_i)\}=\{(t,\psi_i)\}$, it follows readily
from~\eqref{cc89} that the new 1-cocycle $\{\tilde c_{ij} = \tilde b_{ij}/h_i\}$ is cohomologous to the 1-cocycle
$\{tc_{ij}\}$.
Thus the assignment $\{(a_{ij},b_{ij})\} \mapsto \{c_{ij}=b_{ij}/h_i\}$ induces a~map $\rho_{\mathscr U}^{-1}(\xi)\to
H^1(\mathscr U,\mathscr O(D))/\mathbb{C}^*$, where $\mathbb{C}^*$ acts on $H^1(\mathscr U,\mathscr O(D))$ as the scalar
multiplication.
Conversely the assignment $\{c_{ij}\}\mapsto \{(a_{ij},h_ic_{ij})\}$ induces a~map $H^1(\mathscr U,\mathscr O(D))\to
\rho_{\mathscr U}^{-1}(\xi)$, which descends (by~\eqref{cc89}) to a~map from $H^1(\mathscr U,\mathscr
O(D))/\mathbb{C}^*$.
The last map is clearly the inverse of the above map $\rho_{\mathscr U}^{-1}(\xi)\to H^1(\mathscr U,\mathscr O(D))$.
Thus, under f\/ixing $\{a_{ij}\}$ for $\xi\in H^1(\mathscr U,\mathscr O^*)$ and $\{h_i\}$ satisfying $h_i = a_{ij}h_j$,
we obtained a~bijection
\begin{gather}
\label{af2.5}
\rho_{\mathscr U}^{-1}(\xi)
\; \stackrel{\sim}{\longrightarrow}
\;
H^1(\mathscr U,\mathscr O(D))/\mathbb{C}^*.
\end{gather}
Here the quotient space $H^1(\mathscr U,\mathscr O(D))/\mathbb{C}^*$ is of course a~single point if $H^1(\mathscr
U,\mathscr O(D))=0$ and otherwise a~single point plus a~projective space.
The single point corresponds to the line bundle~$\xi$ itself.
In~\cite{At55} it was proved that this map is independent of the choice of $\{a_{ij}\}$ and $\{h_i\}$.
So by taking the inductive limit in~\eqref{af2.5} with respect to open coverings, we obtain a~bijection
\begin{gather*}
\rho^{-1}(\xi)
\;\stackrel{\sim}{\longrightarrow}
\;
H^1(X,\mathscr O(D))/\mathbb{C}^*.
\end{gather*}
Thus from~\eqref{star1} there is a~natural 1-1 correspondence
\begin{gather}
\label{af3}
H^1(X,\mathscr Af)
\;\stackrel{\sim}{\longrightarrow}
\;
\bigsqcup_{\xi\in H^1(X,\mathscr O^*)} H^1(X, \mathscr O(\xi)) /\mathbb{C}^*.
\end{gather}

When $X=\mathbb{CP}^1$, by natural isomorphisms $H^1(X,\mathscr O^*) \simeq H^2(\mathbb{CP}^1,\mathbb{Z})\simeq
\mathbb{Z}$,~\eqref{af3} can be rewritten as
\begin{gather}
\label{af4}
H^1(\mathbb{CP}^1,\mathscr Af)
\;\stackrel{\sim}{\longrightarrow}
\;
\bigsqcup_{n\in\mathbb{Z}} H^1(\mathbb{CP}^1, \mathscr O(n)) /\mathbb{C}^*.
\end{gather}

\begin{Definition}\label{def:degree}
The {\em degree} of an af\/f\/ine bundle $A\to\mathbb{CP}^1$ is the image of~$A\in H^1(\mathbb{CP}^1,\mathscr Af)$ by
the composition of the natural map and identif\/ications
\begin{gather*}
H^1(\mathbb{CP}^1,\mathscr Af)
\;\stackrel{\rho}{\longrightarrow}
\;
H^1(\mathbb{CP}^1,\mathscr O^*)\simeq
H^2(\mathbb{CP}^1,\mathbb{Z})\simeq \mathbb{Z}.
\end{gather*}
\end{Definition}

Evidently, if an af\/f\/ine bundle is a~line bundle, its degree coincides with the usual degree as a~line bundle.
Denoting $\overline A$ for the $\mathbb{CP}^1$-bundle naturally associated to~$A$ as before, the degree of~$A$ is
exactly the negative of the self-intersection number of $\overline A\backslash A$ in $\overline A$.

If $n\ge -1$ we have $H^1(\mathbb{CP}^1,\mathscr O(n))=0$ and $\rho^{-1}(n)$ consists of a~single point which is exactly
the line bundle $\mathscr O(n)$.
Hence if the degree of a~line bundle is more than $-2$, it cannot be deformed even as an af\/f\/ine bundle.
So we are mainly interested in the case where the degree is less than $-1$.
We write such line bundles in the form $\mathscr O(-n)$, so that $n\ge 2$.
Then as $h^1(\mathscr O(-n)) = n-1>0$, $\rho^{-1}(-n)$ consists of the single point (which corresponds to the line
bundle $\mathscr O(-n)$) and the projective space $\mathbb{CP}^{n-2}$ (which is also a~point if $n=2$).
For each $n\ge 2$ we now compute transition law for f\/iber coordinates on arbitrary af\/f\/ine bundles with degree $-n$
in a~concrete form.
For this we take the standard covering $\mathscr U_0:=\{U_0,U_1\}$ where $U_0 =\{(z:w)\in\mathbb{CP}^1\,|\, z\neq 0\}$
and $U_1 =\{(z:w)\in\mathbb{CP}^1\,|\, w\neq 0\}$, and put $u = w/z$, $v=1/u$.
Then we have $H^1(\mathscr U_0,\mathscr O^*)\simeq H^1(\mathbb{CP}^1, \mathscr O^*)$, and so for $\xi=\mathscr O(-n)$,
as $\{a_{ij}\}$ and $\{h_i\}$ we can take
\begin{gather*}
a_{01} = \frac 1{u^n},
\qquad
h_0 =1
\qquad
(\text{so that}
\quad
a_{10} = u^n,
\quad
h_1 = {u^n}).
\end{gather*}
Moreover as a~basis of $H^1(\mathscr U_0,\mathscr O(D)) \!=\!H^1(\mathscr U_0,\mathscr O(-n))\!\simeq\! \mathbb{C}^{n-1}$ we can
take $\{u^{-1},u^{-2},\dots,u^{1-n}\}$ where $u^{-l}\in H^0(U_{01},\mathscr O(D))$.
Thus for each element $b_{01}=t_1u^{-1}+t_2u^{-2} + \dots + t_{n-1} u^{1-n}\in H^1(\mathscr U_0,\mathscr O(D))$ we can
associate an af\/f\/ine $\mathbb{C}$-bundle $A\to\mathbb{CP}^1$ whose transition law for f\/iber coordinates is given by
\begin{gather}
\label{af5}
\zeta_0 = \frac 1{u^n} \zeta_1 + \sum\limits_{l=1}^{n-1} \frac {t_l} {u^l}
\qquad
\text{on}
\quad
U_{01},
\end{gather}
where $\zeta_0$ and $\zeta_1$ are f\/iber coordinates over $U_0=\mathbb{C}(u)$ and $U_1=\mathbb{C}(v)$ respectively.
This can be regarded as def\/ining a~holomorphic family of af\/f\/ine $\mathbb{C}$-bundles over $\mathbb{CP}^1$
parametrized by $\mathbb{C}^{n-1}$, and over the origin we have the line bundle $\mathscr O(-n)$.
We write the total space of this family by~$\mathscr A_n$, thereby obtaining a~holomorphic map
\begin{gather}
\label{af6}
\mathscr A_n\to\mathbb{C}^{n-1}.
\end{gather}
The equation~\eqref{af5} is linear in the variables $\zeta_0$, $\zeta_1$, $t_1,t_2,\dots,t_{n-1}$.
Hence the total space of the family~\eqref{af6} admits a~$\mathbb{C}^*$-action which is the multiplication to all these
variables by weight~$1$.
This $\mathbb{C}^*$-action clearly descends to the scalar multiplication on the base space $\mathbb{C}^{n-1}$, and
f\/ibers of the family~\eqref{af6} are mutually isomorphic along orbits of this $\mathbb{C}^*$-action on
$\mathbb{C}^{n-1}$.
This explains geometrically why two cocycles $\{(a_{ij},b_{ij})\}$ and $\{(a_{ij},tb_{ij})\}$ with values in $\mathscr
Af$ determines mutually isomorphic af\/f\/ine bundles.

As above, for any $t\in\mathbb{C}^*$, the two cocycles $\{(a_{ij}, b_{ij})\}$ and $\{(a_{ij}, tb_{ij})\}$ determines the
same (or isomorphic, more precisely) af\/f\/ine bundles.
This can also be seen directly by noticing that the equation~\eqref{af5} is linear in the variables
$\zeta_0$, $\zeta_1$, $t_1,t_2,\dots,t_{n-1}$.
Hence by identifying f\/ibers of~\eqref{af6} lying over the same linear 1-dimensional subspace, we have obtained
a~family of af\/f\/ine $\mathbb{C}$-bundles over $\mathbb{CP}^1$ which is parametrized by $H^1(\mathbb{CP}^1,\mathscr
O(-n))/\mathbb{C}^*$.
By varying $n$ in $\mathbb{Z}$, this gives a~concrete realization of the bijection~\eqref{af4}.

Strictly speaking, in the argument of the last paragraph, we need to show that the natural map $H^1(\mathscr
U_0,\mathscr Af)\to H^1(\mathbb{CP}^1,\mathscr Af)$ is bijective; especially we need to show that any af\/f\/ine bundle
over $\mathbb{CP}^1$ can be trivialized over the open sets $U_0$ and $U_1$ respecting the structure of af\/f\/ine
bundle.
But this can be proved by standard adjusting argument using coboundaries, and we omit the detail.

\subsection[Af\/f\/ine $\mathbb{C}$-bundles over $\mathbb{CP}^1$ and Hirzebruch surfaces]{Af\/f\/ine
$\boldsymbol{\mathbb{C}}$-bundles over $\boldsymbol{\mathbb{CP}^1}$ and Hirzebruch surfaces}
\label{ss:AH}

We are concerned with ALE SFK metrics on the total spaces of af\/f\/ine $\mathbb{C}$-bundles over $\mathbb{CP}^1$ whose
degree is negative.
We will investigate this through the natural compactif\/ication of the af\/f\/ine bundles to $\mathbb{CP}^1$ bundles.
The latter are of course Hirzebruch surfaces.
In this subsection we will brief\/ly explain relationship between af\/f\/ine $\mathbb{C}$-bundles over $\mathbb{CP}^1$
and the Hirzebruch surfaces.

First let $A\to \mathbb{CP}^1$ be an af\/f\/ine $\mathbb{C}$-bundle, and let $\overline A\to \mathbb{CP}^1$ be the
natural compactif\/ication to a~$\mathbb{CP}^1$-bundle induced by the inclusion ${\rm Af}(\mathbb{C})\subset
\mathrm{PGL}(2,\mathbb{C})$ as before.
We have $\overline A\simeq \mathbb{F}_n$ for some $n\ge 0$.
We write $L:= \overline A\backslash A$ for the added locus, which is of course a~section of the projection $\overline A\to \mathbb{CP}^1$.
Then $-L^2$ is exactly the degree of $A\to \mathbb{CP}^1$.
Hence any af\/f\/ine $\mathbb{C}$-bundle over $\mathbb{CP}^1$ of negative degree is naturally identif\/ied with the
complement of a~section of some $\mathbb{F}_n$ whose self-intersection number is positive.
We write by
\begin{gather}
\label{Kf}
\mathscr F_n\to \mathbb{C}^{n-1}
\end{gather}
for the family of Hirzebruch surfaces that is obtained as the simultaneous compactif\/ication for members of the family
$\mathscr A_n\to \mathbb{C}^{n-1}$ in~\eqref{af6}.
For this family, it is well-known that the Kodaira--Spencer map
\begin{gather}
\label{ks001}
T_0\mathbb{C}^{n-1}\to H^1(\mathbb{F}_n,\Theta)
\end{gather}
at the origin is isomorphic (see~\cite[pp.~309--312]{K}), and the family~\eqref{Kf} gives the Kuranishi family of the
Hirzebruch surface $\mathbb{F}_n$.
Thus the parameter space $\mathbb{C}^{n-1}$ of $\mathscr A_n\to\mathbb{C}^{n-1}$ and $\mathscr F_n\to\mathbb{C}^{n-1}$
may also be naturally identif\/ied with~$H^1(\Theta_{\mathbb{F}_n})$.

Although the transition law for each member of the family~\eqref{Kf} is concretely given as in~\eqref{af5}, it is not
easy to identify them with $\mathbb{F}_m$ for a~precise value of $m$.
This was intensively studied in~\cite[p.~143, Theorem]{Suwa}, where an explicit answer was given, but it is too
complicated to write the result here.
Some exceptions are identif\/ication for f\/ibers on the coordinate axes of~$\mathbb{C}^{n-1}$.
Namely letting $\mathbb{C}(t_l)$ to be the $l$-th coordinate axis of~$\mathbb{C}^{n-1}$, if we introduce new f\/iber
coordinates $\tilde{\zeta}_0$ and $\tilde{\zeta}_1$ by
\begin{gather}
\label{zeta'}
\tilde{\zeta}_0 = \frac{u^l\zeta_0-t_l}{t_l\zeta_0}
\qquad
\text{and}
\qquad
\tilde{\zeta}_1= \frac{\zeta_1}{t_lv^{n-l}\zeta_1 + t_l^2}
\end{gather}
on the open sets $U_0=\mathbb{C}(u)$ and $U_1=\mathbb{C}(v)$ respectively, then with the aid of~\eqref{af5}, we readily
obtain the relation $\tilde{\zeta}_0 = v^{n-2l}\tilde{\zeta}_1$, which means that the ruled surface over the axis
$\mathbb{C}(t_l)$ is isomorphic to $\mathbb{F}_{n-2l}$, except the central f\/iber.
Here we are allowing the case $n-2l<0$ and in that case $\mathbb{F}_{n-2l}$ means $\mathbb{F}_{2l-n}$.
We also note that, as a~natural extension of the $\mathbb{C}^*$-action on $\mathscr A_n$, the total space of $\mathscr
F_n\to\mathbb{C}^{n-1}$ has a~$\mathbb{C}^*$-action, and it also descends to the scalar multiplication on $\mathbb{C}^{n-1}$.

Conversely if~$L$ is a~section of $\pi:\mathbb{F}_m\to\mathbb{CP}^1$ for some $m\ge 0$, then the complement
$\mathbb{F}_m\backslash L$ is biholomorphic to an af\/f\/ine $\mathbb{C}$-bundle over $\mathbb{CP}^1$.
This can be seen in the following way.
Let~$\Gamma_0$ and~$\Gamma_{\infty}$ be sections satisfying $\Gamma_0^2 = -m$ and $\Gamma_{\infty}^2 = m$.
If $m=0$, we assume $\Gamma_0\neq \Gamma_{\infty}$.
Let $\mathscr U=\{U_i\}$ be an open covering of $\mathbb{CP}^1$ which satisf\/ies for any~$i$ at least one of $L\cap
\Gamma_0\cap \pi^{-1} (U_i)=\varnothing$ or $L\cap \Gamma_{\infty}\cap \pi^{-1} (U_i)=\varnothing$ holds for any~$i$.
Let $\zeta_i$ be any f\/iber coordinate over $U_i$ of the line bundle $\mathscr O(-m)\subset\mathbb{F}_m$ (so that
$\Gamma_0\cap\pi^{-1} (U_i)$ and $\Gamma_{\infty}\cap\pi^{-1} (U_i)$ are def\/ined by $\zeta_i=0$ and $\zeta_i=\infty$
respectively), and $f_i$ be a~meromorphic function on $U_i$ such that $L\cap\pi^{-1} (U_i)$ is def\/ined by $\zeta_i =
f_i$.
From the choice, $f_i$ does not have both a~zero and a~pole.
Then for any~$i$ such that $f_i$ does not have a~pole, we def\/ine a~new f\/iber coordinate $\tilde{\zeta}_i$ over $U_i$
as an af\/f\/ine bundle by setting
\begin{gather}
\label{af8}
\tilde{\zeta}_i = \frac 1{\zeta_i-f_i}.
\end{gather}
Then from the choice of~$i$, this may be used as a~f\/iber coordinate on $\mathbb{F}_m\to\mathbb{CP}^1$, and we have
$L\cap \pi^{-1}(U_i) = \{\tilde{\zeta}_i = \infty\}$.
For the remaining~$i$-s, $f_i$ does not have a~zero.
$L\cap \Gamma_0 \cap \pi^{-1} (U_i)=\varnothing$.
For these~$i$-s we put
\begin{gather}
\label{af9}
\tilde{\zeta}_i = \frac{f_i\zeta_i}{f_i-\zeta_i}.
\end{gather}
Then this can also be used as a~f\/iber coordinate over $U_i$, and we again have $L\cap \pi^{-1}(U_i) =
\{\tilde{\zeta}_i = \infty\}$.
From~\eqref{af8} and~\eqref{af9} we readily see that the transition law for the new coordinate system
$\{\tilde{\zeta}_i\}$ is included in the af\/f\/ine group ${\rm Af}(\mathbb{C})$.
Therefore $\mathbb{F}_m\backslash L$ is actually an af\/f\/ine bundle.
However, even if the equation for a~section~$L$ is given in a~concrete form, it is not immediate again to trivialize the
af\/f\/ine bundle $\mathbb{F}_m\backslash L\to\mathbb{CP}^1$ over $U_0=\mathbb{C}(u)$ and $U_1=\mathbb{C}(v)$ and write
down the transition function in the form~\eqref{af5}.

\subsection{Computations for Hirzebruch surfaces}

As we mentioned we will investigate ALE SFK metrics on the af\/f\/ine $\mathbb{C}$-bundles over $\mathbb{CP}^1$ through
the compactif\/ication to Hirzebruch surfaces.
More precisely the Hirzebruch surfaces are included in the twistor spaces of a~conformal compactif\/ication of the ALE
SFK metrics on the af\/f\/ine bundles, and the added section will be the twistor line over the added point at
inf\/inity, whose self-intersection number in the surface is positive.
For this purpose, in this subsection, we make computations for pairs $(\mathbb{F}_n,L)$ where~$L$ is a~section
satisfying $L^2>0$.
(So $\mathbb{F}_n\backslash L$ is an af\/f\/ine $\mathbb{C}$-bundle over $\mathbb{CP}^1$ of negative degree.) Especially
we compute the dimension $h^i(\Theta_{\mathbb{F}_n,L})$ for arbitrary pairs, where $\Theta_{\mathbb{F}_n,L}$ is the
sheaf of germs of holomorphic vector f\/ields on $\mathbb{F}_n$ which are tangent to~$L$.

If~$L$ is a~section of $\mathbb{F}_n\to\mathbb{CP}^1$ which satisf\/ies $L^2>0$, we have $L\in |\Gamma_0 + (n+l)f|$ for
some $l\ge 0$ (see Section~\ref{ss:notation} for notation). Of course the value of $h^i(\Theta_{\mathbb{F}_n,L})$
depends on the number $l$.
We begin with the case $l=0$.
In this case we have $n>0$ as we are supposing $L^2>0$.
Moreover if we identify $\mathbb{F}_n$ with $\overline A$ where $A=\mathscr O(-n)$, then~$L$ can be identif\/ied with
$\overline A\backslash A$, the section at inf\/inity.
Hence regarding $\mathscr O(-n)$ (or $\mathbb{F}_n$) as the minimal resolution of $\mathbb{C}^2/\mathbb{Z}_n$ (or
$\mathbb{CP}^2/\mathbb{Z}_n$) where $\mathbb{Z}_n\subset{\rm{GL}}(2,\mathbb{C})$ is a~cyclic subgroup of scalar matrices
of order $n$, the pair $(\mathbb{F}_n,L)$ has an ef\/fective action of ${\rm{GL}}(2,\mathbb{C})/\mathbb{Z}_n$.
In particular we have $h^0(\Theta_{\mathbb{F}_n,L}) \ge 4$.
\begin{Proposition}
\label{prop:l0}
Suppose $n>0$ and let~$L$ be any $(+n)$-section of $\mathbb{F}_n\to\mathbb{CP}^1$.
Then we have the following:
\begin{enumerate}\itemsep=0pt
\item[$(i)$] $h^0(\Theta_{\mathbb{F}_n, L})=4$ and $H^2(\Theta_{\mathbb{F}_n, L})=0$.
\item[$(ii)$] The natural homomorphism $H^1(\Theta_{\mathbb{F}_n, L}) \to H^1(\Theta_{\mathbb{F}_n})$ is isomorphic,
and these are $(n-1)$-dimensional vector spaces.
\item[$(iii)$] The complex structure of the pair $(\mathbb{F}_n,L)$ is independent of the choice of~$L$.
\item[$(iv)$] We have ${\rm{Aut}}_0(\mathbb{F}_n,L)\simeq {\rm{GL}}(2,\mathbb{C})/\mathbb{Z}_n$ $($see above$)$.
\end{enumerate}
\end{Proposition}

\begin{proof}
Let $\Theta_{\mathbb{F}_n/\mathbb{CP}^1}\subset\Theta_{\mathbb{F}_n}$ be the subsheaf consisting of germs of holomorphic
vector f\/ields which are tangent to f\/ibers of $\pi:\mathbb{F}_n\to\mathbb{CP}^1$.
Scalar matrices in ${\rm{GL}}(2,\mathbb{C})$ induce a~$\mathbb{C}^*$-action on $\mathbb{F}_n$ which preserves each
f\/iber of $\mathbb{F}_n\to\mathbb{CP}^1$, and it def\/ines a~vector f\/ield which is tangent to each f\/iber of $\pi$.
Hence we obtain a~section of $\Theta_{\mathbb{F}_n/\mathbb{CP}^1}$.
Moreover, as the vector f\/ield has simple zeros on $\Gamma_0\sqcup L$ and no other zeros, we obtain
$\Theta_{\mathbb{F}_n/\mathbb{CP}^1}\simeq \mathscr O_{\mathbb{F}_n}(\Gamma_0+ L )$.
With the aid of this isomorphism, we have the following commutative diagram of exact sequences of sheaves on
$\mathbb{F}_n$:
\begin{gather*}
\begin{CD}
& & 0 & & 0 & & 0 & &
\\
& & @VVV @VVV @VVV
\\
0 @>>> \mathscr O_{\mathbb{F}_n} (\Gamma_0) @>>> \Theta_{\mathbb{F}_n/\mathbb{CP}^1} @>>> N_{ L /\mathbb{F}_n} @>>> 0
\\
& & @VVV @VVV @|
\\
0 @>>> \Theta_{\mathbb{F}_n, L } @>>> \Theta_{\mathbb{F}_n} @>>> N_{ L /\mathbb{F}_n} @>>> 0
\\
& & @VVV @VVV @VVV
\\
0 @>>> \pi^*\Theta_{\mathbb{CP}^1} @= \pi^*\Theta_{\mathbb{CP}^1} @>>> 0
\\
& & @VVV @VVV
\\
&& 0 && 0
\end{CD}
\end{gather*}
Since the group ${\rm{GL}}(2,\mathbb{C})$ acts transitively on $\mathbb{CP}^1$, the natural map
$H^0(\Theta_{\mathbb{F}_n,L})\to H^0(\pi^*\Theta_{\mathbb{CP}^1})$ is surjective.
Hence from the f\/irst column of the diagram, since $h^0(\mathscr O_{\mathbb{F}_n}(\Gamma_0)) = 1$ as $n>0$ and
$h^0(\Theta_{\mathbb{CP}^1}) = 3$, we obtain $h^0(\Theta_{\mathbb{F}_n,L}) = 1 + 3 = 4$.
Also from the same column, as we readily have $h^1(\mathscr O_{\mathbb{F}_n}(\Gamma_0))=n-1$ and $h^2(\mathscr
O_{\mathbb{F}_n}(\Gamma_0))=0$, we obtain $h^1(\Theta_{\mathbb{F}_n,L}) = n-1$ and $h^2(\Theta_{\mathbb{F}_n,L}) = 0$.

From the isomorphism $\Theta_{\mathbb{F}_n/\mathbb{CP}^1}\simeq \mathscr O_{\mathbb{F}_n}(\Gamma_0+ L )$ we also obtain
$h^0(\Theta_{\mathbb{F}_n/\mathbb{CP}^1}) = h^0(\mathscr O_{\mathbb{F}_n}(\Gamma_0+ L )) = n + 2$.
Hence as $h^0(N_{ L /\mathbb{F}_n}) =h^0(\mathscr O_{\mathbb{CP}^1}(n)) = n+1$ and $h^0(\mathscr
O_{\mathbb{F}_n}(\Gamma_0)) = 1$, from the f\/irst row, we obtain that the map $H^0( \Theta_{\mathbb{F}_n/\mathbb{CP}^1}
) \to H^0(N_{ L /\mathbb{F}_n})$ is surjective.
Hence from the commutative diagram the map $H^0( \Theta_{\mathbb{F}_n} ) \to H^0(N_{ L /\mathbb{F}_n})$ is also
surjective.
Therefore from the middle row, as $h^0(\Theta_{\mathbb{F}_n,L}) = 4 $, we obtain $h^0(\Theta_{\mathbb{F}_n}) = 4 + (n+1)
= n+5$.
Moreover, as $h^i(N_{L/\mathbb{F}_n}) = h^i(\mathscr O(n))=0$ for $i\in\{1,2\}$ and $h^2(\Theta_{\mathbb{F}_n,L}) = 0$,
we obtain from the same row that the natural map $H^1(\Theta_{\mathbb{F}_n,L})\to H^1(\Theta_{\mathbb{F}_n})$ is
isomorphic and $h^2(\Theta_{\mathbb{F}_n})=0$.

The assertion (iii)
is clear since $\mathbb{F}_n\backslash L$ is isomorphic to the line bundle $\mathscr O(-n)$ for any
$(+n)$-section~$L$.
(iv) follows from the remark preceding to Proposition~\ref{prop:l0} and the assertion (i) which is already shown.
\end{proof}

Next we consider the case $l=1$, which requires some care.

\begin{Proposition}\label{prop:l1}
Suppose $n>0$.
If~$L$ is a~section of $\pi:\mathbb{F}_n\to\mathbb{CP}^1$ belonging to the system $|\Gamma_0 + (n+1)f|$, we have the following:
\begin{enumerate}\itemsep=0pt
\item[$(i)$] $ h^0(\Theta_{\mathbb{F}_n,L})=2$, $h^1(\Theta_{\mathbb{F}_n,L})= n-1$, and $h^2(\Theta_{\mathbb{F}_n,L})= 0$.
\item[$(ii)$] The natural map $H^1(\Theta_{\mathbb{F}_n,L})\to H^1(\Theta_{\mathbb{F}_n})$ is isomorphic.
\item[$(iii)$] The complex structure of the pair $(\mathbb{F}_n,L)$ is independent of the choice of the section~$L$.
\end{enumerate}
\end{Proposition}

\begin{proof}
Let $L\in |\Gamma_0 +(n+1)f|$ be a~section as in the proposition.
Then we readily have
\begin{gather}
\label{inter1n1}
(L,\Gamma_0)=1,
\qquad
(L,\Gamma_0+nf) = n+1.
\end{gather}
Note that we have not specif\/ied a~$(+n)$-section $\Gamma_{\infty}$ yet.
We write $p=L\cap\Gamma_0$ and let $q\in L$ be any point which is dif\/ferent from $p$.
Then as $\dim |\Gamma_0+nf| = n+1$, by dimension counting, there exists a~section $\Gamma_{\infty} \in |\Gamma_0 + nf|$
which touches~$L$ at the point $q$ by multiplicity $(n+1)$.
Let $T_{\mathbb{C}}\subset{\rm{Aut}}\mathbb{F}_n$ be the maximal torus which is determined by the property that it
preserves the two sections $\Gamma_0$, $\Gamma_{\infty}$ and f\/ixes the two points $p$, $q$.
The complement $\mathbb{F}_n\backslash \Gamma_{\infty}$ may be identif\/ied with the line bundle $\mathscr O(-n)$.
Let $u$ be an af\/f\/ine coordinate on $\mathbb{CP}^1\backslash \pi(q)$ (where $\pi$ is the projection
$\mathbb{F}_n\to\mathbb{CP}^1$ as before), and $\zeta$ a~f\/iber coordinate of the line bundle $\mathscr O(-n)$ over
$U_0=\mathbb{C}(u)$, so that
\begin{gather*}
p=(0,0),
\qquad
q=(\infty,\infty),
\qquad
\Gamma_{0} = \{\zeta=0\},
\qquad
\text{and}
\qquad
\Gamma_{\infty} = \{\zeta=\infty\}.
\end{gather*}
Then as~$L$ intersects $\Gamma_0$ transversally at $p$ by~\eqref{inter1n1}, a~def\/ining equation for~$L$ has to be of
the form, in the above coordinates,
\begin{gather*}
\zeta = uh(u),
\qquad
h(0)\neq 0,
\end{gather*}
where $h=h(u)$ is a~holomorphic function on $U_0 = \mathbb{C}(u)$.
In the coordinates $(v,\eta^{-1}):= (u^{-1},u^{-n}\zeta^{-1})$ around the point $q=(\infty,\infty)$, this can be
rewritten as $ \eta^{-1} = v^{n+1}/h(v^{-1})$.
Then since~$L$ touches $\Gamma_{\infty}$ at the point $q$ by multiplicity $(n+1)$, the function $h(v^{-1})$ cannot have
a~pole at $v=0$.
This means that $h(u)$ is a~constant.
Hence~$L$ is def\/ined by the equation $\zeta = cu$ for some $c\in \mathbb{C}^*$.
But we may assume $c=1$ by changing the f\/iber coordinate $\zeta$ to $c^{-1}\zeta$.
Thus~$L$ is def\/ined by $\zeta=u$ in the coordinates $(u,\zeta)$.
In particular this proves the assertion (iii).

Next in order to determine $h^0(\Theta_{\mathbb{F}_n,L})$, we recall that, in terms of the coordinates $(u,\zeta)$, any
vector f\/ield $\theta\in H^0(\Theta_{\mathbb{F}_n})$ is concretely written as (see~\cite[pp. 43--44]{MKbook})
\begin{gather}
\label{MK1}
\theta= g(u) \frac{\partial}{\partial u} + \left( f(u) \zeta^ 2 + c\zeta u \right) \frac{\partial}{\partial \zeta},
\end{gather}
where $g(u) = a_1 u^2 + a_2 z + a_3$
$(a_i\in\mathbb{C})$,
$f(u) = b_1 u^n + b_2 u^{n-1} + \dots + b_{n+1}$
$(b_i\in\mathbb{C})$ and $c\in \mathbb{C}$.
(So we have $3 + (n+1) + 1 = n + 5$ parameters in total, which agrees with $h^0(\Theta_{\mathbb{F}_n}) = n+5$.)

For later use, we let $l\ge 1$ and let the section~$L$ to be def\/ined by $\zeta = u^l$, and def\/ine $F(u,\zeta):=
\zeta - u^l$, so that $F$ is a~def\/ining equation of~$L$.
Then $\theta\in H^0(\Theta_{\mathbb{F}_n,L})$ if\/f the derivation $\theta F$ satisf\/ies $\theta F|_L = 0$.
By~\eqref{MK1} we have
\begin{gather*}
\theta F = f\zeta^2 + c u \zeta - lgu^{l-1}.
\end{gather*}
Hence by substituting $\zeta = u^l$, the restriction becomes
\begin{gather}
\theta F |_L = \big(b_1 u^n + b_2 u^{n-1} + \dots + b_{n+1}\big) u^{2l} + c u^{l+1} - l \big(a_1 u^2 + a_2 u + a_3\big) u^{l-1}
\nonumber
\\
\phantom{\theta F |_L}{}
= \big(b_1 u^{2l+n} + b_2 u^{2l+n-1} + \dots + b_{n+1}u^{2l}\big) + \big\{(c-la_1) u^{l+1} - la_2 u^l - la_3 u^{l-1}\big\}.
\label{der1}
\end{gather}
When $l=1$, we have $2l = l+1$, and we obtain
\begin{gather*}
\theta F|_L = \big(b_1 u^{n+2} + b_2 u^{n+1} + \dots + b_{n+1}u^{2}\big) + \big\{(c-a_1) u^{2} - a_2 u - a_3 \big\}
\\
\phantom{\theta F|_L}{}
= b_1 u^{n+2} + b_2 u^{n+1} + \dots + b_n u^3 + \{b_{n+1} + (c -a_1)\} u^2 - a_2 u - a_3.
\end{gather*}
Thus $\theta F|_L = 0$ if\/f{\samepage
\begin{gather*}
b_1=b_2=\dots=b_{n} = 0,
\qquad
b_{n+1} + c - a_1 = a_2 = a_3 = 0.
\end{gather*}
These imply $h^0(\Theta_{\mathbb{F}_n, L}) = 2$.}

It remains to compute $h^i(\Theta_{\mathbb{F}_n, L})$ for $i\in \{1,2\}$ and show the isomorphicity in~(ii).
But these follow readily from~\eqref{cohohir}, the standard exact sequence $ 0 \to \Theta_{\mathbb{F}_n, L} \to
\Theta_{\mathbb{F}_n} \to N_{L/\mathbb{F}_n}\to 0$ and $h^0(\Theta_{\mathbb{F}_n, L})=2$.
The assertion (iii) is already shown.
\end{proof}

By the proposition, the group ${\rm{Aut}}_0(\mathbb{F}_n,L)$ is 2-dimensional when $l=1$.
This group can also be readily determined in a~concrete form.
For this, as before let $p$ be the intersection point of~$L$ and $\Gamma_0$.
(By~\eqref{inter1n1}~$L$ and $\Gamma_0$ intersect transversally at a~unique point.)

\begin{Proposition}\label{prop:nonred}
Suppose $n>0$ and let~$L$ be any section of $\pi:\mathbb{F}_n\to \mathbb{CP}^1$ belonging to the system $|\Gamma_0+(n+1)f|$.
Then the $2$-dimensional group ${\rm{Aut}}_0(\mathbb{F}_n, L)$ can be naturally identified with
the group $\{g\in\mathrm{PGL}(2,\mathbb{C}) \,|\, g(\pi(p)) = \pi(p)\}$, which is isomorphic
to the affine transformation group ${\rm Af}(\mathbb{C})$.
\end{Proposition}

\begin{proof}
The kernel sheaf of the restriction of the natural surjection $\Theta_{\mathbb{F}_n}\to \pi^*\Theta_{\mathbb{CP}^1}$ to
the subsheaf~$\Theta_{\mathbb{F}_n,L}$ can be obtained in a~similar way to the f\/irst column of the commutative diagram
in the proof of Proposition~\ref{prop:l0}, and consequently we obtain the exact sequence
\begin{gather}
\label{ntses}
0
 \longrightarrow
\mathscr O_{\mathbb{F}_n} (\Gamma_0 - f)
 \longrightarrow
\Theta_{\mathbb{F}_n,L}
 \longrightarrow
\pi^*\Theta_{\mathbb{CP}^1}
\longrightarrow
0.
\end{gather}
Clearly we have $H^0 (\mathscr O_{\mathbb{F}_n} (\Gamma_0 - f ) )=0$.
Hence we obtain that the natural homomorphism $H^0(\Theta_{\mathbb{F}_n,L} ) \to H^0( \pi^*\Theta_{\mathbb{CP}^1})\simeq
H^0(\Theta_{\mathbb{CP}^1})$ is injective.
Therefore, unlike ${\rm{Aut}}_0\mathbb{F}_n$, the subgroup ${\rm{Aut}}_0(\mathbb{F}_n,L)$ can be regarded as
a~subgroup of $\mathrm{PGL}(2,\mathbb{C})$.
Moreover any $g\in {\rm{Aut}}_0(\mathbb{F}_n,L)$ has to f\/ix the point $p = \Gamma_0\cap L$, since $\Gamma_0$ is
${\rm{Aut}}\mathbb{F}_n$-invariant as $n>0$.
Hence under the above inclusion ${\rm{Aut}}_0(\mathbb{F}_n,L)\subset\mathrm{PGL}(2,\mathbb{C})$,
${\rm{Aut}}_0(\mathbb{F}_n,L)$ is included in the subgroup of $\mathrm{PGL}(2,\mathbb{C})$ in the proposition.
But since we already know $\dim {\rm{Aut}}_0(\mathbb{F}_n,L)=2$ by Proposition~\ref{prop:l1}, we obtain the coincidence.
\end{proof}

Thus the computations for $h^i(\Theta_{\mathbb{F}_n,L})$ and ${\rm{Aut}}_0(\mathbb{F}_n,L)$ is over for arbitrary
sections when $l\in \{0,1\}$.
Next we consider the case $l>1$.
In this case, the situation is not completely homogeneous:

\begin{Proposition}\label{prop:l2}
Suppose $n>0$, $l>1$, and let~$L$ be any section of $\mathbb{F}_n\to \mathbb{CP}^1$ belonging to the system $|\Gamma_0+(n+l)f|$.
Then we have $H^2(\Theta_{\mathbb{F}_n,L})=0$.
Further one of the following holds:
\begin{enumerate}\itemsep=0pt
\item[$(i)$] $h^0(\Theta_{\mathbb{F}_n, L}) = 1$, $h^1(\Theta_{\mathbb{F}_n, L}) = n+2l-4$, and there is an exact
sequence
\begin{gather*}
0
\longrightarrow
\mathbb{C}^{2l-3}
\longrightarrow
H^1(\Theta_{\mathbb{F}_n,L})
\longrightarrow
H^1(\Theta_{\mathbb{F}_n})
\longrightarrow
0.
\end{gather*}
\item[$(ii)$] $h^0(\Theta_{\mathbb{F}_n, L}) = 0$, $h^1(\Theta_{\mathbb{F}_n, L}) = n+2l-5$, and there is an exact
sequence
\begin{gather*}
0
\longrightarrow
\mathbb{C}^{2l-4}
\longrightarrow
H^1(\Theta_{\mathbb{F}_n,L})
\longrightarrow
H^1(\Theta_{\mathbb{F}_n})
\longrightarrow
0.
\end{gather*}
\end{enumerate}
Furthermore, as long as the section~$L$ satisfies~$(i)$, the complex structure of the pair $(\mathbb{F}_n,L)$ is
independent of the choice of~$L$, and we have ${\rm{Aut}}_0(\mathbb{F}_n, L) \simeq \mathbb{C}^*$.
\end{Proposition}

\begin{proof}
The vanishing $H^2 (\Theta_{\mathbb{F}_n, L}) =0$ is an immediate consequence of the exact sequence
\begin{gather}
\label{ses:56}
0
\longrightarrow
\Theta_{\mathbb{F}_n,L}
\longrightarrow
\Theta_{\mathbb{F}_n}
\longrightarrow
N_{L/\mathbb{F}_n}
\longrightarrow
0
\end{gather}
since as $N_{L/\mathbb{F}_n}\simeq\mathscr O(n+2l)$ we have $H^1(N_{L/\mathbb{F}_n}) = 0$ and also
$H^2(\Theta_{\mathbb{F}_n}) =0$.

{\sloppy The ingredient is to show $h^0(\Theta_{\mathbb{F}_n, L}) \le 1$.
For this we f\/irst note that, in the same way to the f\/irst column of the commutative diagram in the proof of
Proposition~\ref{prop:l0} or the exact sequence~\eqref{ntses} in the case $l=1$, we have an exact sequence
\begin{gather*}
0
\longrightarrow
\mathscr O_{\mathbb{F}_n} (\Gamma_0 - lf )
\longrightarrow
\Theta_{\mathbb{F}_n,L}
\longrightarrow
\pi^*\Theta_{\mathbb{CP}^1}
\longrightarrow
0.
\end{gather*}
This again means that the natural map $H^0(\Theta_{\mathbb{F}_n,L} ) \to H^0(\Theta_{\mathbb{CP}^1})$ is injective, and
hence ${\rm{Aut}}_0(\mathbb{F}_n,L)$ may be considered as a~subgroup of $\mathrm{PGL}(2,\mathbb{C})$.
Moreover, since $\Gamma_0$ is ${\rm{Aut}}\mathbb{F}_n$-invariant, ele\-ments of ${\rm{Aut}}_0(\mathbb{F}_n,L)$ f\/ix any
point of the intersection $\Gamma_0\cap L$.
Since we have $(L,\Gamma_0)=(\Gamma_0 + (n+l)f,\Gamma_0) = l\ge 2$, $L\cap \Gamma_0$ is non-empty.
If it consists of more than two points, then the image of ${\rm{Aut}}_0(\mathbb{F}_n,L)\to \mathrm{PGL}(2,\mathbb{C})$
is clearly identity, and so ${\rm{Aut}}_0(\mathbb{F}_n,L)$ is trivial, meaning $h^0(\Theta_{\mathbb{F}_n,L})=0$.

}

If $L\cap \Gamma_0$ consists of two points, the image of ${\rm{Aut}}_0(\mathbb{F}_n,L)\to \mathrm{PGL}(2,\mathbb{C})$ is
included in the $\mathbb{C}^*$-subgroup determined by the two points.
Suppose that the image is actually the $\mathbb{C}^*$-subgroup, and let $T_{\mathbb{C}}$ be the maximal torus of
${\rm{Aut}}_0(\mathbb{F}_n)$ which contains ${\rm{Aut}}_0(\mathbb{F}_n,L)(\simeq \mathbb{C}^*)$.
Then $T_{\mathbb{C}}$ determines on $\mathbb{F}_n$ a~structure of toric surface, and singles out a~$(+n)$-section
$\Gamma_{\infty}$ by $T_{\mathbb{C}}$-invariance.
Moreover~$L$ cannot intersect $\Gamma_{\infty}$ since~$L$ minus the two f\/ixed points $L\cap\Gamma_0$ forms an orbit of
the $\mathbb{C}^*$-subgroup of $T_{\mathbb{C}}$, and $\Gamma_{\infty}$ is disjoint from the unique 2-dimensional orbit
of the $T_{\mathbb{C}}$-action.
This contradicts $(L,\Gamma_{\infty})=n+l$ $(>0)$.
Therefore if $L\cap \Gamma_0$ consists of two points, ${\rm{Aut}}_0(\mathbb{F}_n,L)$ is trivial.

If $L\cap \Gamma_0$ consists of one point, since $(L,\Gamma_0)=l$, as in the same way to the proof of
Proposition~\ref{prop:l1}, we can f\/ind coordinates $(u,\zeta)$ on the line bundle $\mathscr O(-n)\subset\mathbb{F}_n$
such that the point $L\cap \Gamma_0$ corresponds to the origin and~$L$ is def\/ined by an equation $\zeta = u^l$.
In these coordinates the $(+n)$-section def\/ined by the equation $\zeta=\infty$ intersects~$L$ at the unique point
$(\infty,\infty)$ by the biggest multiplicity $(n+l)$.
Then using the computations in Proposition~\ref{prop:l1}, by writing a~vector f\/ield $\theta\in H^0(\Theta_{\mathbb{F}_n})$
as in~\eqref{MK1}, we have~\eqref{der1} for the derivative $\theta F|_L$ of the def\/ining equation $F=\zeta-u^l$ of~$L$.
Now as $l>1$ we have $2l>l+1$.
Hence looking the powers to $u$ in~\eqref{der1}, the vanishing $\theta F|_L = 0$ is equivalent to the equations
\begin{gather*}
b_1=b_2=\dots=b_{n+1} = 0,
\qquad
c-la_1 = a_2 = a_3 = 0.
\end{gather*}
From these we obtain $h^0(\Theta_{\mathbb{F}_n, L})=1$.
Thus we have seen that if $l>1$ we always have \mbox{$h^0(\Theta_{\mathbb{F}_n,L}) \le 1$} and the equality holds exactly
when~$L$ touches the section $\Gamma_0$ at a~point by the biggest multiplicity.

Once this is obtained, the assertions (i) and (ii) are readily obtained from the exact sequence~\eqref{ses:56}.
We omit the detail.
The f\/inal assertion is clear from the above argument since we have seen that if $h^0(\Theta_{\mathbb{F}_n,L}) = 1$,
equation for~$L$ can be taken as $\zeta=u^l$ in the coordinates $(u,\zeta)$ on the line bundle $\mathscr O(-n)$.
\end{proof}

As a~corollary to the results in this subsection, we obtain the following result on the existence and uniqueness up to
isomorphisms of $\mathbb{C}^*$-invariant pair $(\mathbb{F}_n,L)$ for each $n$ and $l$:

\begin{Corollary}\label{cor:unique}
For each integers $n\ge 0$ and $l\ge 0$, there exists a~section~$L$ of $\mathbb{F}_n\to\mathbb{CP}^1$ which satisfies
the following two properties:
\begin{enumerate}\itemsep=0pt
\item[$(i)$] $L\in |\Gamma_0+(n+l)f|$,
\item[$(ii)$] $L$ is $\mathbb{C}^*$-invariant, where $\mathbb{C}^*$ is a~subgroup of ${\rm{Aut}}\mathbb{F}_n$ which acts
non-trivially on~$L$.
\end{enumerate}
Moreover, for each $n$ and $l$ the complex structure of the pair $(\mathbb{F}_n,L)$ is independent of the choice of such
a~section~$L$.
\end{Corollary}

\begin{proof}
The assertion for the case $n>0$ follows from Propositions~\ref{prop:l0},~\ref{prop:l1} and~\ref{prop:l2}.
The assertion for the case $n=0$ is immediate to see.
\end{proof}

Finally in this subsection, we discuss variation of the complex structures of the af\/f\/ine bundles in the family
$\mathscr A_n\to\mathbb{C}^{n-1}$.
As in the beginning of this subsection, by identifying the total space of the line bundle $\mathscr O(-n)$ with the
minimal resolution of $\mathbb{C}^2/\mathbb{Z}_n$, the line bundle $\mathscr O(-n)$ admits
a~${\rm{GL}}(2,\mathbb{C})$-action.
This naturally gives rise to a~${\rm{GL}}(2,\mathbb{C})$-action on the cohomology group $H^1(\mathbb{CP}^1,\mathscr
O(-n))$.
Recalling that this cohomology group is exactly the base space of the family $\mathscr A_n\to\mathbb{C}^{n-1}$, f\/ibers
of the family are mutually biholomorphic if they are over the same orbit of the ${\rm{GL}}(2,\mathbb{C})$-action.

From the results in~\cite{HonCMP2}, the ${\rm{GL}}(2,\mathbb{C})$-action on the base space $\mathbb{C}^{n-1}$ is
identif\/ied with the tensor product
\begin{gather*}
S_1^{n-2}\mathbb{C}^2:=S^{n-2}\mathbb{C}^2\otimes\mathbb{C}_1,
\end{gather*}
where $S^{n-2}\mathbb{C}^2$ is the $(n-2)$-th symmetric product of the natural ${\rm{GL}}(2,\mathbb{C})$-action on
$\mathbb{C}^2$, and $\mathbb{C}_1$ is the 1-dimensional representation of ${\rm{GL}}(2,\mathbb{C})$ which is just the
multiplication of the determinant.
(See Section~\ref{explicit}, especially the isomorphisms~\eqref{isom49}.) It follows that if $n\in\{2,3\}$ the
${\rm{GL}}(2,\mathbb{C})$-action on the base space $\mathbb{C}^{n-1}\backslash\{0\}$ is transitive.
Therefore when $n\in\{2,3\}$, any member of the family $\mathscr A_n\to\mathbb{C}^{n-1}$ is biholomorphic except the
central f\/iber $\mathscr O(-n)$.
These can also be seen by just noting that, if $n\in\{2,3\}$, any f\/iber except the central f\/iber is identif\/ied
with $\mathbb{F}_{n-2}\backslash L$ with $L\in |\Gamma_0+(n-1)f|$, and also that the complex structure of the pair
$(\mathbb{F}_{n-2},L)$ is independent of the choice of a~non-singular member $L\in|\Gamma_0+(n-1)f|$ by
Proposition~\ref{prop:l1}.

On the other hand, if $n\ge 4$, the ${\rm{GL}}(2,\mathbb{C})$-action on the base space
$S_1^{n-2}\mathbb{C}^2=\mathbb{C}^{n-1}$ minus the origin is not transitive, and so the quotient space
$(\mathbb{C}^{n-1}\backslash\{0\})/{\rm{GL}}(2,\mathbb{C})$ consists of more than two elements.
As the total spaces of the af\/f\/ine bundles are not only open but also do not have compact holomorphic curves, it
seems dif\/f\/icult to determine when two af\/f\/ine surfaces lying over dif\/ferent ${\rm{GL}}(2,\mathbb{C})$-orbits
are mutually biholomorphic (if $n\ge 4$).

\subsection{Computations for surfaces of smooth normal crossing}
\label{ss:sncs}

In this subsection we f\/irst construct a~variety of smooth normal crossing from two copies of the Hirzebruch surface
$\mathbb{F}_n$ by identifying the same sections, and then compute cohomology groups for them.
In the next section these varieties will be included as a~subvariety in twistor spaces of the 4-dimensional orbifold
$\widehat{\mathscr O(-n)}$.

For this let $D$ be a~non-singular complex surface and $L\subset D$ a~non-singular rational curve.
Denoting $J$ for the complex structure on~$D$, we denote by $\overline D$ the complex surface obtained from~$D$ by
changing the complex structure~$J$ to $-J$.
Let $\id :D\to \overline D$ be the identity map.
This is an {\em anti-}holomorphic map.
Write $\overline L:=\id(L)\subset \overline D$.
Let $\tau:L\to L$ be an anti-holomorphic involutions of~$L\simeq\mathbb{CP}^1$, and we def\/ine a~map
$\phi:L\to\overline L$ by $\phi:=\id|_{L}\circ\tau$.
Since both~$\id|_{L}$ and~$\tau$ are anti-holomorphic, $\phi$ is a~holomorphic map.
Let $D\cup_{L,\tau}\overline D$ be the space obtained from the disjoint union $D\sqcup \overline D$ by identifying~$L$
and~$\overline L$ by~$\phi$.
By the holomorphicity of~$\phi$, $D\cup_{L,\tau}\overline D$ is naturally equipped with the structure of a~complex
variety which is smooth normal crossing.

Let $\sigma:D\sqcup \overline D\to D\sqcup \overline D$ be the map def\/ined by
\begin{gather*}
\sigma (p) =
\begin{cases}
\id (p) &\text{if}\quad p\in D,
\\
\id^{-1}(p) &\text{if}\quad p\in \overline D.
\end{cases}
\end{gather*}
This is clearly an involution which f\/lips the two components (as the map $\id$ f\/lips from the def\/inition),
and is an anti-holomorphic map since $\id$ and $\id^{-1}$ are.
If the two points $p\in L$ and $q\in \overline L$ satisfy $q = \phi(p)$, we have
\begin{gather*}
\phi(\sigma(q)) = \phi(\id^{-1}(q))=\id\circ\tau\circ\id^{-1}(q)=\id\circ\tau\circ\id^{-1} (\id\circ\tau(p)) = \id(p)=\sigma(p).
\end{gather*}
Namely we have $\phi(\sigma(q)) = \sigma(p)$.
Hence $\sigma$ descends to an endomorphism of $D\cup_{L,\tau}\overline D$.
We use the same letter $\sigma$ for this map.
This is an anti-holomorphic involution since the original $\sigma$ is.
Thus the variety $D\cup_{L,\tau}\overline D$ is naturally equipped with a~real structure.
The structure of $D\cup_{L,\tau}\overline D$ as a~complex variety with a~real structure depends not only on the
rational curve~$L$ but also on the involution $\tau$.
Further, if $p\in L$, we have
\begin{gather*}
\phi^{-1}(\sigma(p)) = (\id\circ\tau)^{-1}(\id (p)) = \tau^{-1}(p) = \tau (p).
\end{gather*}
This means that on the intersection $D\cap \overline D\subset D\cup_{L,\tau}\overline D$, the involution $\sigma$ may be
identif\/ied with the involution $\tau$ on~$L$.

We apply this construction to the pair $(D,L) = (\mathbb{F}_n,L)$, where $n > 0$ and~$L$ is a~section of
$\mathbb{F}_n\to \mathbb{CP}^1$ with a~positive self-intersection number and an anti-holomorphic involution $\tau:L\to L$ without a~f\/ixed point.
As above the structure of the resulting variety $\mathbb{F}_n\cup_{L,\tau}\overline{\mathbb{F}}_n$ depends on the choice
of the involution $\tau$.
But if the section~$L$ is supposed to be invariant under a~$\mathbb{C}^*$-action on $\mathbb{F}_n$ that acts on
non-trivially on~$L$, then the choice of $\tau$ is naturally constrained to be $\mathbb{C}^*$-equivariant, and
consequently if $p$ and $q$ denote the f\/ixed points of the $\mathbb{C}^*$-action on~$L$, we have $\tau(p)\in\{p,q\}$.
But since $\tau$ is supposed to have no f\/ixed point, we obtain $\tau(p) = q$.
This means that in an af\/f\/ine coordinate $u$ on~$L$ for which the $\mathbb{C}^*$-action is given by $u\mapsto tu$ for
$t\in\mathbb{C}^*$, we can write $\tau(u) = -a/\overline u$ for some $a>0$.
Therefore the ef\/fect of varying $\tau$ (namely varying the number $a>0$) is absorbed in the $\mathbb{C}^*$-action on
$\mathbb{F}_n$, and moreover by Corollary~\ref{cor:unique}, the complex structure of the pair $(\mathbb{F}_n,L)$ is
independent of the choice of such a~section~$L$.
Consequently the variety $\mathbb{F}_n\cup_{L,\tau}\overline{\mathbb{F}}_n$ makes a~unique sense.
Further the $\mathbb{C}^*$-actions on $(\mathbb{F}_n,L)$ and $(\overline{\mathbb{F}}_n,\overline L)$ are naturally glued
and the variety is equipped with a~$\mathbb{C}^*$-action.
As we are particularly interested in these varieties, we introduce notation for them:

\begin{Definition}
\label{def:scn}
For integers $n \ge 0$ and $l\ge 0$, let $L\in |\Gamma_0 + (n+l)f|$ be any $\mathbb{C}^*$-invariant section on
$\mathbb{F}_n$, and we denote by $\mathbb{F}_n\cup_l\overline{\mathbb{F}}_n$ for the variety of simple normal crossing
with $\mathbb{C}^*$-action, which is obtained from the two copies of the pair $(\mathbb{F}_n,L)$ by identifying two
$L$-s by an anti-holomorphic involution $\tau$ without a~f\/ixed point in the above way.
\end{Definition}

The notation $\mathbb{F}_n\cup_l\overline{\mathbb{F}}_n$ ref\/lects the independency from the choices of~$L$ and $\tau$.
Thus the complex structure of this variety is solely determined by two non-negative integers $n$ and $l$.
For these varieties we have the following.

\begin{Proposition}\label{prop:fd2}
Let $n>0$ and $l\ge 0$.
Then for the tangent sheaf $\Theta$ of the variety $\mathbb{F}_n \cup_l \mathbb{F}_n$ above, we have the following:
\begin{itemize}\itemsep=0pt
\item[$(i)$] If $l=0$, we have
\begin{gather*}
h^0(\mathbb{F}_n \cup_l \overline{\mathbb{F}}_n,\Theta) = 5,
\qquad
h^1(\mathbb{F}_n \cup_l \overline{\mathbb{F}}_n,\Theta) = 2(n-1),
\qquad
h^2(\mathbb{F}_n \cup_l \overline{\mathbb{F}}_n,\Theta) = 0.
\end{gather*}
\item[$(ii)$] If $l\ge 1$, we have
\begin{gather*}
h^0(\mathbb{F}_n \cup_l \overline{\mathbb{F}}_n,\Theta) = 1,
\qquad
h^1(\mathbb{F}_n \cup_l \overline{\mathbb{F}}_n,\Theta) = 2(n+2l-3),
\qquad
h^2(\mathbb{F}_n \cup_l \overline{\mathbb{F}}_n,\Theta) = 0.
\end{gather*}
\end{itemize}
\end{Proposition}

\begin{proof}
Though these can be shown in a~standard way by using Propositions~\ref{prop:l0},~\ref{prop:l1} and~\ref{prop:l2}, we
write a~proof as there is a~subtle point that relies on our construction of the variety
$\mathbb{F}_n\cup_l\overline{\mathbb{F}}_n$.
We have the standard exact sequence $0 \to \Theta_{\mathbb{F}_n\cup_l \overline{\mathbb{F}}_n} \to
\Theta_{\mathbb{F}_n,L}\oplus\Theta_{\overline{\mathbb{F}}_n,\overline L}\to \Theta_L\to 0$, where $L\in
|\Gamma_0+(n+l)f|$ is a~$\mathbb{C}^*$-invariant section identif\/ied by $\phi$.
For the case $l=0$, the natural map $H^0(\Theta_{\mathbb{F}_n,L}) \to H^0(\Theta_{L})$ is surjective from that of
${\rm{Aut}}_0(\mathbb{F}_n,L)\to \Aut L$.
Therefore from the above exact sequence and Proposition~\ref{prop:l0} (i), (ii) we obtain the required value for $h^i
(\mathbb{F}_n \cup_0 \overline{\mathbb{F}}_n,\Theta)$ for any~$i$ as well as natural isomorphisms {\samepage
\begin{gather}
H^1(\mathbb{F}_n \cup_0 \overline{\mathbb{F}}_n,\Theta)  \simeq H^1(\mathbb{F}_n, \Theta_{\mathbb{F}_n,L}) \oplus H^1
(\overline{\mathbb{F}}_n,\Theta_{\overline{\mathbb{F}}_n,\overline L})
\label{isom001}
\\
\hphantom{H^1(\mathbb{F}_n \cup_0 \overline{\mathbb{F}}_n,\Theta)}{}
\simeq H^1(\mathbb{F}_n, \Theta_{\mathbb{F}_n}) \oplus H^1 (\overline{\mathbb{F}}_n,\Theta_{\overline{\mathbb{F}}_n}).
\label{isom002}
\end{gather}}

\noindent
Next  for the case $l=1$, the natural homomorphism ${\rm{Aut}}_0(\mathbb{F}_n,L)\to \Aut L$ is not surjective and the
image is the af\/f\/ine transformation group as in Proposition~\ref{prop:nonred}.
Namely it consists of elements of $\Aut L$ which f\/ixes the point $p=\Gamma_0\cap L$.
As our involution $\tau$ is supposed to interchange the two f\/ixed points $p$ and $q$ of the $\mathbb{C}^*$-action, it
follows that the image of the natural map $H^0(\mathbb{F}_n, \Theta_{\mathbb{F}_n,L}) \oplus H^0
(\overline{\mathbb{F}}_n,\Theta_{\overline{\mathbb{F}}_n,\overline L})\to H^0(\Theta_L)$ is again surjective since the
two af\/f\/ine groups generate $\Aut L$.
Hence the cohomology exact sequence takes the same form as the case $l=0$, and by using Proposition~\ref{prop:l1} (i)
and (ii), we obtain the required value of $h^i (\mathbb{F}_n \cup_1 \overline{\mathbb{F}}_n,\Theta)$ as well as the
natural isomorphisms~\eqref{isom001} and~\eqref{isom002}.

Finally if $l>1$, by Proposition~\ref{prop:l2} the image of the natural injection ${\rm{Aut}}_0(\mathbb{F}_n,L)\to \Aut
L$ is the $\mathbb{C}^*$-subgroup that f\/ixes the two points $p$ and $q$.
Therefore from our choice of $\tau$, the image of the natural map $H^0(\mathbb{F}_n, \Theta_{\mathbb{F}_n,L}) \oplus H^0
(\overline{\mathbb{F}}_n,\Theta_{\overline{\mathbb{F}}_n,\overline L})\to H^0(\Theta_L)$ is 1-dimensional.
Therefore from the cohomology sequence we obtain $H^0(\Theta_{\mathbb{F}_n\cup_l\overline{\mathbb{F}}_n})
\simeq\mathbb{C}$, the exact sequence
\begin{gather*}
0
\longrightarrow
\mathbb{C}^2
\longrightarrow
H^1 (\mathbb{F}_n \cup_l \overline{\mathbb{F}}_n,\Theta)
\longrightarrow
H^1(\mathbb{F}_n, \Theta_{\mathbb{F}_n,L}) \oplus H^1
(\overline{\mathbb{F}}_n,\Theta_{\overline{\mathbb{F}}_n,\overline L})
\longrightarrow
0
\end{gather*}
and the isomorphism $H^2 (\mathbb{F}_n \cup_l \overline{\mathbb{F}}_n,\Theta) \simeq H^2(\mathbb{F}_n,
\Theta_{\mathbb{F}_n,L}) \oplus H^2 (\overline{\mathbb{F}}_n,\Theta_{\overline{\mathbb{F}}_n,\overline L})$.
From Proposition~\ref{prop:l2}, we f\/inish the proof of the assertion (ii).
\end{proof}

\section{Computations for twistor spaces}
\label{explicit}

In this section, based on the results in the previous section, we intensively study small deformations of the LeBrun
metric on $\mathscr O(-n)$ which preserve ALE SFK properties, and in particular show that any af\/f\/ine
$\mathbb{C}$-bundle over $\mathbb{CP}^1$ of negative degree admits an ALE SFK metric.
Next we investigate small deformations of the metrics on the af\/f\/ine bundles again as ALE SFK metrics, and in
particular show that even if we f\/ix the complex structure on the af\/f\/ine bundles, they admit a~1-parameter
deformation for which the conformal classes are not constant.

\subsection{Some generalities on twistor spaces of ALE SFK metrics} 

Before starting actual computations, we brief\/ly recall basic properties of the twistor spaces of ALE SFK metrics,
including its natural compactif\/ication.
These will be used for investigating deformations of metrics which preserve ALE SFK property.
For more precise treatment on compactif\/ications of ALE ASD 4-manifolds, we refer the paper~\cite{V}.

Let $(X,J)$ be a~complex surface, $g$ an ASD Hermitian metric on it, $Z$ the twistor space of the ASD conformal class
$[g]$, and $F:=K_Z^{-1/2}$ the natural square root of the anticanonical line bundle of $Z$, which is available on any
twistor space.
Then the complex structure $J$ determines a~section of the twistor projection $Z\to X$ in a~tautological way, and its
image becomes a~non-singular divisor $D$ on $Z$.
$D$ is biholomorphic to $X$ by the projection $Z\to X$.
Let $\overline D$ be the divisor determined by the conjugate complex structure $-J$ on $X$.
We always have $D\cap \overline D=\varnothing$ as $J\neq -J$.
Then Pontecorvo's theorem~\cite{Pont92} means that the ASD Hermitian metric $g$ is K\"ahler with respect to $J$ if and
only if $D + \overline D\in |F|$.

If $X$ is non-compact with one end and the ASD metric $g$ is asymptotically Euclidean at inf\/inity, then $(X,[g])$ can
be compactif\/ied as an ASD manifold by adding a~point at inf\/inity.
Let $(\hat X,[\hat g])$ be the resulting compact ASD manifold, and $\hat Z$ the twistor space of $(\hat X, [\hat g])$,
which is smooth.
Then the closure ${\rm{Cl}}(D)$ of the above divisor $D\subset Z$ is a~divisor in $\hat Z$, and from the ALE SFK
property of the metric~\cite[proof of Proposition~6, p.~312]{LB92}, the divisor ${\rm{Cl}}(D)$ satisf\/ies the following
properties: (i)~${\rm{Cl}}(D)$ is still non-singular and ${\rm{Cl}}(D) = D\sqcup L$, where~$L$ is the twistor line over
the point at inf\/inity, (ii)~${\rm{Cl}}(D)\cap {\rm{Cl}}(\overline D)=L$, and the intersection is transverse, and (iii)~the normal bundle of~$L$ in ${\rm{Cl}}(D)$ (and also in ${\rm{Cl}}(\overline D)$) is of degree one.
Conversely a~divisor in $\hat Z$ satisfying these properties determines, up to overall constants, an ASD K\"ahler metric
on $(X,J)$ which is asymptotically Euclidean at inf\/inity.

When the SFK surface $(X,J,g)$ is ALE in a~strict sense (i.e.\ asymptotic to the f\/lat Euclidean orbifold
$\mathbb{C}^2/\Gamma$ at inf\/inity, where $\Gamma$ is a~non-trivial f\/inite subgroup of ${\rm{U}}(2)$ acting freely on
the unit sphere), the pair $(X,g)$ has a~natural compactif\/icaton $(\hat X, \hat g)$ as an ASD orbifold, which means
that~$\hat X$ is an orbifold of the form $X\cup\{\infty\}$ with $\infty$ being an orbifold point of~$\hat X$, and~$\hat
g$ is an ASD orbifold metric on~$\hat X$ whose conformal class on $X$ remains to be equal to~$g$.
Also the twistor space~$Z$ of~$(X,g)$ has a~natural compactif\/ication, for which we again denote by~$\hat Z$.
This is of course the twistor space of the ASD orbifold~$(\hat X,\hat g)$ in a~natural sense, and we again have~$\hat Z
= Z\sqcup L$, where~$L$ is the twistor line over the orbifold point~$\infty$.
We have $\Sing\hat Z\subset L$, and all singularities are quotient singularity by the group which is
orientation-reversing conjugate (namely conjugate after reversing the orientation; see~\cite[Def\/inition~1.4]{V} for the precise
def\/inition) to the above group~$\Gamma$.
Especially, denoting ${\rm{U}}(1)\subset{\rm{U}}(2)$ for the subgroup of consisting of scalar matrices, if~$\Gamma$ is
a~cyclic subgroup of ${\rm{U}}(1)$ with order $n\ge 2$, then $\hat Z$ has $A_{n-1}$-singularities along~$L$.
(This is particular to these subgroups, and for other subgroup $\Gamma\subset U(2)$, singular points of~$\hat Z$ are
isolated.) Moreover if $D\subset Z$ again denotes the divisor determined by the complex structure~$J$ on~$X$ and
${\rm{Cl}}(D)$ means its closure in $\hat Z$, then ${\rm{Cl}}(D)$ itself (and therefore ${\rm{Cl}}(\overline D)$ also)
is a~non-singular (but non-Cartier) divisor on $\hat Z$.
Moreover we have ${\rm{Cl}}(D)\cap {\rm{Cl}}(\overline D)=L$, and the normal bundle of~$L$ in ${\rm{Cl}}(D)$ (and also
in ${\rm{Cl}}(\overline D)$) is of degree $n$.
Furthermore the union ${\rm{Cl}}(D) \cup {\rm{Cl}}(\overline D)$ itself is smooth normal crossing.
We also note that in this situation the natural extension of the line bundle $F$ over $Z$ to $\hat Z$ is not just an
orbifold bundle but an ordinary line bundle; in other words the sum ${\rm{Cl}}(D) + {\rm{Cl}}(\overline D)$ is a~Cartier
divisor on $\hat Z$, while ${\rm{Cl}}(D)$ and ${\rm{Cl}}(\overline D)$ are not.

\looseness=-1
Because $\mathscr O(-n)$ is obtained as the minimal resolution of the quotient space $\mathbb{C}^2/\Gamma$ where
$\Gamma\subset{\rm{U}}(2)$ is the cyclic subgroup of scalar matrices of order $n$, ALE SFK metrics on $\mathscr O(-n)$
give rise to the last situation where the compactif\/ied twistor space $\hat Z$ has $A_{n-1}$-singularities along the
twistor line~$L$ at inf\/inity.
Here, we do not suppose that the complex structure on $\mathscr O(-n)$ is the natural one and we will also consider
complex structures which support the af\/f\/ine bundles in Section~\ref{ss:ab}.
Let $\widehat{\mathscr O(-n)}$ be the one-point compactif\/ication of the 4-manifold $\mathscr O(-n)$, and in the
following, instead of the letters $\hat Z$ and ${\rm{Cl}}(D)$, we use the letters $Z$ and $D$ respectively to mean the
twistor space of the conformal compactif\/ication of an ALE SFK metric on the 4-mani\-fold~$\mathscr O(-n)$ and the
(non-Cartier) divisor on $Z$ determined by the complex structure on the 4-mani\-fold~$\mathscr O(-n)$.
In this situation $D$ is biholomorphic to $\mathbb{F}_{n-2k}$ for some $k\ge 0$ satisfying $n-2k\ge 0$.
This is because $D$ contains the twistor line~$L$ at inf\/inity as a~$(+n)$-curve as above, which means the rationality
of $D$; further the decomposition $D=\mathscr O(-n)\sqcup L$ as a~smooth manifold means $b_2(D) = 2$, and hence
$D\simeq\mathbb{F}_m$ for some $m\ge 0$; but $\mathbb{F}_m\to\mathbb{CP}^1$ has a~$(+n)$-section if\/f $m=n-2k$ for some $k\ge0$.
Of course we have $k=0$ if the complex structure on $\mathscr O(-n)$ is the natural one.
We also remark that if $\tau$ is an anti-holomorphic involution of~$L$ without a~f\/ixed point, the union
$D\cup\overline D$, which is a~Cartier divisor in $Z$ as above, is isomorphic to the surface $D\cup_{L,\tau}\overline D$
constructed in the f\/irst half of Section~\ref{ss:sncs}, as a~complex variety with a~real structure.

Thus if $(Z,D)$ is a~pair of a~compact but singular twistor space and a~divisor determined by an ALE SFK metric on the
4-manifold $\mathscr O(-n)$, deformations of the metric preserving ALE SFK property are equivalent to locally trivial
deformations of the pair $(Z,D\cup\overline D)$ preserving the real structure.
For details on locally trivial deformations for complex spaces and pairs of a~complex space and a~complex subspace of
it, we refer a~book~\cite[Section~3.4]{SeBook}.
In particular, if we def\/ine the subsheaf $\Theta_{Z, D\cup\overline D}$ of the tangent sheaf $\Theta_Z$ by
\begin{gather*}
\Theta_{Z, D\cup\overline D}:=\big\{v\in \Theta_Z\,|\, v(f) \in \mathscr I_{D\cup\overline D}
\text{ if }
f\in \mathscr I_{D\cup\overline D}\big\},
\end{gather*}
then f\/irst order deformations of the pair $(Z,D\cup\overline D)$ which are locally trivial are in one to one
correspondence with the cohomology group $H^1(\Theta_{Z,D\cup\overline D})$, and obstructions are in
$H^2(\Theta_{Z,D\cup\overline D})$.
In particular if $H^2(\Theta_{Z,D\cup\overline D})=0$, the Kuranishi family for locally trivial deformations of the pair
$(Z, D\cup\overline D)$ is constructed over a~neighborhood of the origin in $H^1(\Theta_{Z,D\cup\overline D})$.

\subsection{Deformations of the LeBrun metric}
\label{ss:DLm}

Having recalled these basic materials, we start to investigate deformations of the LeBrun's ALE SFK metric on $\mathscr
O(-n)$ as an ALE SFK metrics, by investigating locally trivial deformations of the pair $(Z,D\cup\overline D)$ of
compactif\/ied singular twistor space and the divisor.
The following proposition provides basic information about such deformations.
\begin{Proposition}\label{prop:isom1}
Suppose $n\ge 3$ and let $Z$ be the twistor space on the orbifold $\widehat{\mathscr O(-n)}$, which is associated to the
conformal compactification of the LeBrun's ALE-SFK metric on $\mathscr O(-n)$ with negative mass.
Let $D$ be the divisor on $Z$ which is the closure of the section of the twistor fibration that is determined by the
complex structure of $\mathscr O(-n)$.
$(D$ is biholomorphic to $\mathbb{F}_n.)$ Then we have
\begin{gather}
\label{cohomdim1}
H^i( \Theta_Z(-D-\overline D) ) = 0,
\quad
i\neq 1,
\qquad
H^1( \Theta_Z(-D-\overline D) ) \simeq \mathbb{C},
\\
H^2(\Theta_{Z, D\cup\overline D}) = H^2 ( \Theta_{D\cup \overline D} ) = 0.
\label{van1}
\end{gather}
Moreover there is a~natural isomorphism
\begin{gather}
\label{basicisom1}
H^1 (\Theta_{Z,  D \cup \overline D} ) \simeq H^1 (D, \Theta_{ D \cup \overline D} ),
\end{gather}
and these are $2(n-1)$-dimensional.
Furthermore the natural map
\begin{gather}
\label{forget}
H^1(\Theta_{Z,D\cup\overline D})
\longrightarrow
H^1(\Theta_Z)
\end{gather}
is injective, and if $n=3$, this is moreover surjective.
\end{Proposition}

The isomorphism~\eqref{basicisom1} will be of fundamental importance in the rest of this article.

\begin{proof}
The vanishing $H^2( \Theta_Z(-D-\overline D) ) = 0$ immediately follows from~\cite[Proposition~3.1]{HonCMP2} since $S$ in the
proposition is a~divisor in the system $|F|$ and hence $\Theta_Z(-D-\overline D)\simeq\Theta_Z(-S)$.
In order to compute $h^i(\Theta_Z(-D-\overline D) )$ for $i\in \{0,1,3\}$,
we use computations in the proof of the above proposition in~\cite{HonCMP2}.
Noting $\Theta_Z (-D-\overline D ) \simeq \Theta_Z\otimes F^{-1}$, the isomorphisms~(3.5),~(3.6) and~(3.8) in the proof
of~\cite[Proposition~3.1]{HonCMP2} are valid not only for $H^2$ but also for $H^i$ for any~$i$ because we have $H^i(\mathscr
O_{\mathbb{CP}^1} (-1))=H^i({\mathbb{CP}^1\times\mathbb{CP}^1},\mathscr O(-1,-1))=0$ for any~$i$.
Therefore in the notation of that proof, we have $H^i(\Theta_Z\otimes F^{-1}) \simeq H^i(X,\mathscr L')$ for any~$i$.
Further from the exact sequence~(3.9) there, we have $H^i(\mathscr L')\simeq H^i(\mathscr F')$ for any~$i$.
Furthermore from the exact sequence~(3.10) there, we obtain
\begin{gather*}
H^i(\mathscr F') = 0,
\quad
i\neq 1,
\qquad
H^1(\mathscr F')\simeq H^0(\Delta,\mathscr O)  (\simeq\mathbb{C}).
\end{gather*}
These in particular imply~\eqref{cohomdim1}.

Next in order to deduce~\eqref{van1} and~\eqref{basicisom1} we consider the standard exact sequence
\begin{gather}
\label{ses:1}
0
\longrightarrow
\Theta_Z(-D-\overline D)
\longrightarrow
\Theta_{Z, D\cup\overline D}
\longrightarrow
\Theta_{D\cup \overline D}
\longrightarrow
0.
\end{gather}
Since the isometry group of the LeBrun metric on $\mathscr O(-n)$ is $U(2)/\mathbb{Z}_n$ \cite{LB88}, where
$\mathbb{Z}_n$ is the cyclic subgroup consisting of scalar matrices of order $n$, we have $ h^0 ( \Theta_{Z,
D\cup\overline D}) = 4.
$ On the other hand, as $D\cup \overline D\simeq \mathbb{F}_n\cup_0\overline{\mathbb{F}}_n$ biholomorphically, by
Proposition~\ref{prop:fd2}~(i) we have
\begin{gather}
\label{cohomdim2}
h^0 ( \Theta_{D\cup\overline D}) = 5,
\qquad
h^1 ( \Theta_{D\cup\overline D}) = 2(n-1),
\qquad
h^2 ( \Theta_{D\cup\overline D}) = 0.
\end{gather}
Therefore using~\eqref{cohomdim1} the cohomology exact sequence of~\eqref{ses:1} implies
\begin{gather}
\label{les:1}
0
\longrightarrow
H^0(\Theta_{Z, D\cup\overline D})
\big({\simeq}\,\mathbb{C}^4\big)
\longrightarrow
H^0(\Theta_{D\cup \overline D})
\big({\simeq}\, \mathbb{C}^5\big)
\nonumber
\\
\hphantom{0}{}
\longrightarrow
H^1(\Theta_Z(-D-\overline D))
(\simeq \mathbb{C})
\longrightarrow
H^1(\Theta_{Z, D\cup\overline D})
\longrightarrow
H^1(\Theta_{D\cup \overline D})
\big({\simeq}\, \mathbb{C}^{2(n-1)}\big)
\nonumber
\\
\hphantom{0}{}
\longrightarrow
0
\longrightarrow
H^2(\Theta_{Z, D\cup\overline D})
\longrightarrow
0.
\end{gather}
From this we obtain $H^2( \Theta_{Z, D\cup\overline D}) = 0$, an exact sequence
\begin{gather}
\label{LB:key}
0
\longrightarrow
H^0(\Theta_{Z,D\cup\overline D})
\big({\simeq}\, \mathbb{C}^4\big)
\longrightarrow
H^0(\Theta_{D\cup \overline D})
\big({\simeq}\, \mathbb{C}^5\big)
\longrightarrow
H^1( \Theta_Z(-D-\overline D))
(\simeq \mathbb{C})
\longrightarrow
0,\!\!\!\!
\end{gather}
and also the isomorphism~\eqref{basicisom1}.
From the last isomorphism we obtain $h^1(\Theta_{Z,D\cup\overline D}) = 2(n-1)$ by~\eqref{cohomdim2}.

Finally we show that the map~\eqref{forget} is injective.
For this let $N'$ be the cokernel sheaf of the natural injection $\Theta_{Z,D\cup\overline D}\to \Theta_Z$.
We have an exact sequence $0\to \Theta_{Z, D\cup \overline D}\to \Theta_Z\to N'\to 0$, and so for the injectivity it
suf\/f\/ices to show $H^0(N') = 0$.
Let $N:=\mathscr O_Z(D+\overline D)|_{D\cup\overline D}$ be the normal sheaf of the divisor $D\cup \overline D$ in~$Z$.
Since $D+\overline D\in |F|$ and $F$ is an ordinary line bundle on~$Z$, the sheaf~$N$ is an invertible $\mathscr
O_{D\cup\overline D}$-module, and isomorphic to $F|_{D+\overline D}$.
Then by computing local generators of the sheaves $\Theta_{Z}$ and $\Theta_{Z,D\cup\overline D}$ in coordinates, and
then comparing the resulting generators of the cokernel sheaf~$N'$ with local generators of~$N$, we obtain a~natural
isomorphism
\begin{gather}
\label{nottrivial}
N'\simeq N\otimes_{\mathscr O_{D\cup\overline D}}\mathscr I_{L},
\end{gather}
where $L=D\cap \overline D$ is the twistor line over the point at inf\/inity as before, and $\mathscr I_L$ is the ideal
sheaf of~$L$ in $D\cup \overline D$.
On the other hand, by the adjunction formula we have $K_{D\cup\overline D} \simeq K_Z + [D+\overline D]|_{D\cup\overline
D}\simeq -2F + F|_{D+\overline D}\simeq -F|_{D\cup\overline D}$.
Hence from~\eqref{nottrivial} we obtain $N'\simeq -K_{D\cup\overline D}\otimes \mathscr I_L$.
Further for the canonical sheaf of $D\cup\overline D$, as this itself is smooth normal crossing, we have
\begin{gather*}
K_{D\cup\overline D}|_D \simeq K_D +[\overline D]|_D \simeq K_D+\mathscr O_D(L),
\end{gather*}
and similar for $K_{D\cup\overline D}|_{\overline D}$.
Hence by taking the inverse for these and taking a~tensor product with~$\mathscr O(-L)$, we obtain
\begin{gather}
\label{N''}
N'|_D\simeq -K_D-\mathscr O_D(2L),
\qquad
N'|_{\overline D}\simeq -K_{\overline D}-\mathscr O_{\overline D}(2L).
\end{gather}
Now as $D\simeq\mathbb{F}_n$ we have $-K_D\simeq 2\Gamma_0+(n+2)f$, and as~$L$ is a~$(+n)$-section
we have $\mathscr O_D(L)\simeq \Gamma_0+nf$.
Hence we have
\begin{gather}
\label{N'}
-K_D-\mathscr O_D(2L) \simeq 2\Gamma_0+(n+2)f -2(\Gamma_0 + n f) \simeq -(n-2)f.
\end{gather}
Thus as $n-2>0$ from the assumption $n>2$ we obtain $H^0(-K_D-\mathscr O_D(2L))=0$.
With reality, this means $H^0(D\cup \overline D,N')=0$.
Thus the injectivity of~\eqref{forget} follows.
If $n=3$, the map is also surjective since we have $h^1(\Theta_Z) = 4(n-2) = 4$ by~\cite[Proposition~2.1]{HonCMP2}, which
coincides with $2(n-1) = 4$.
\end{proof}

\begin{Remark}
The computations and the conclusions in the proposition are valid also for the case $n=2$ except the injectivity of the
map~\eqref{forget}.
For the case $n=2$, as in~\eqref{N'}, we have $N'|_D\simeq \mathscr O_D$.
With reality this means $N'\simeq \mathscr O_{D\cup\overline D}$, and hence we have $H^0(N')\simeq\mathbb{C}$.
Further from the cohomology exact sequence this is mapped to $H^1(\Theta_{Z,D\cup \overline D})$ injectively.
Thus the map~\eqref{forget} has a~1-dimensional kernel.
\end{Remark}

Next, letting $Z$ and $D$ be as in Proposition~\ref{prop:isom1}, we collect basic results on versal families of locally
trivial deformations of $Z$, $(Z,D\cup\overline D)$ and $D\cup \overline D$ and their relationship, which are readily
derived from Proposition~\ref{prop:isom1} and the results in Section~\ref{s:Hirz}.

First, for the the twistor space $Z$ of the LeBrun structure on $\widehat{\mathscr O(-n)}$, as showed
in~\cite[Proposition~2.1]{HonCMP2}, we have $H^2(\Theta_Z)=0$ and $h^1(\Theta_Z)=4n-8$.
Hence the parameter space of the Kuranishi family of locally trivial deformations of $Z$ may be identif\/ied with
a~neighborhood of the origin in $H^1(\Theta_Z)\simeq\mathbb{C}^{4n-8}$.
Versal family of~$Z$ as twistor spaces is obtained as the restriction of the Kuranishi family onto the real locus of the
neighborhood.
We denote the last real locus by~$\mathscr U_{\rm{ASD}}$, which is clearly smooth and real $(4n-8)$-dimensional.
As in~\cite{HonCMP2} we call the corresponding family of ASD conformal structures on $\widehat{\mathscr O(-n)}$
(parameterized by $\mathscr U_{\rm{ASD}}$) as the {\em versal family of ASD structures} for the LeBrun structure.
If $n>3$, not all these ASD structures preserve the K\"ahler representative.
From the construction we have a~canonical isomorphism
\begin{gather*}
T_0\mathscr U_{\rm{ASD}}\simeq H^1(\Theta_Z)^{\sigma}
\end{gather*}
as real vector spaces, where the upper-script means the real subspace.

Second, for the pair $(Z,D\cup\overline D)$, in a~similar way to the above argument, since $ H^2 (\Theta_{Z,  D \cup
\overline D} ) =0$ and $ h^1 (\Theta_{Z,  D \cup \overline D} ) =2(n-1)$ as in Proposition~\ref{prop:isom1}, the
parameter space of the Kuranishi family for locally trivial deformations of the pair $(Z,D\cup\overline D)$ is
identif\/ied with a~neighborhood of the origin in $H^1 (\Theta_{Z, D \cup \overline D} ) \simeq\mathbb{C}^{2(n-1)}$.
Restricting this to the real locus, we obtain a~deformation of $Z$ preserving not only a~structure of twistor space but
also the K\"ahler representative in the conformal class.
Let $\mathscr K'\subset H^1(\Theta_{Z,D\cup\overline D})^{\sigma}$ be the parameter space of this family.
We have a~natural isomorphism $ T_0\mathscr K'\simeq H^1(\Theta_{Z,D\cup\overline D})^{\sigma} $.

For a~relationship between the families over $\mathscr U_{\rm{ASD}}$ (of twistor spaces) and $\mathscr K'$ (of pairs of
twistor spaces and Cartier divisors), by versailty, after a~possible shrinking of the domain, there is an induced map,
for which we denote by $\psi_1$, from $\mathscr K'$ to $\mathscr U_{\rm{ASD}}$, such that the pullback by~$\psi_1$ of
the family over $\mathscr U_{\rm{ASD}}$ is isomorphic to the $Z$-portion of the family of pairs over~$\mathscr K'$.
Though $\psi_1$ is not uniquely determined, the derivative $\psi'_1(0)$ is exactly the restriction of the
map~\eqref{forget} to the real locus.
By the proposition the last map is injective, and moreover isomorphism if $n=3$.
So if $n=3$ we may think $\mathscr K' \simeq \mathscr U_{\rm{ASD}}$ by $\psi_1$.
If $n>3$, since $h^1(\Theta_{Z,D\cup\overline D}) = 2(n-1)$ by Proposition~\ref{prop:isom1}, $\psi_1:\mathscr
K'\to\mathscr U_{\rm{ASD}}$ is an embedding as a~real submanifold of dimension $2(n-1)$ in~$\mathscr U_{\rm{ASD}}$ (and
$\dim \mathscr U_{\rm{ASD}} = 4(n-2)$ as above).
We call the image $\psi_1(\mathscr K')$ the {\em K\"ahler locus} in~$\mathscr U_{\rm ASD}$ and denote it by $\mathscr
K$.
If $n=3$, we may think $\mathscr K = \mathscr U_{\rm{ASD}}$ as above.
From the construction we have a~natural isomorphism
\begin{gather*}
\psi_1'(0)^{-1}: \ T_0\mathscr K
\;\stackrel{\sim}{\longrightarrow}
\;
H^1(\Theta_{Z,D\cup\overline D})^{\sigma}
\end{gather*}
as real vector spaces, where the upper-script means the real subspace.

Next for locally trivial deformations of the variety $D\cup\overline D$, since $H^2(\Theta_{D\cup\overline D}) = 0$ by
Proposition~\ref{prop:fd2} (i), the Kuranishi family is parameterized by a~neighborhood of the origin in
$H^1(\Theta_{D\cup\overline D})$.
Denote $\mathscr J\subset H^1(\Theta_{D\cup\overline D})^{\sigma}$ for the real locus of the neighborhood.
Then again by versality, after a~possible shrinking of the domain, there is an induced map, for which we denote by
$\psi_2$, from~$\mathscr K'$ to~$\mathscr J$ that induces an isomorphism between the two families.
Similarly to $\psi_1$, while~$\psi_2$ is not uniquely determined, the derivative $\psi'_2(0)$ is identif\/ied with the
real part of the natural map $H^1(\Theta_{Z,D\cup\overline D})\to H^1(\Theta_{D\cup\overline D})$.
The last map is an isomorphism by~\eqref{basicisom1}, and therefore $\psi_2$ is isomorphic in a~neighborhood of the
origin in~$\mathscr K'$.
Hence the composition $\psi_2\circ\psi_1^{-1}$ gives an isomorphism from the K\"ahler locus $\mathscr K\subset\mathscr
U_{\rm{ASD}}$ to $\mathscr J$, and the $D\cup\overline D$-portion of the families of pairs over~$\mathscr K'$ and the
family over $\mathscr J$ are isomorphic by~$\psi_2$.
The situation is summarized as in the following diagram:
\begin{gather}
\begin{CD}
H^1(\Theta_{D\cup\overline D})^{\sigma} @<{\sim}<< H^1(\Theta_{Z,D\cup\overline D})^{\sigma} @>{\rm{inj.}}>>
H^1(\Theta_{Z})^{\sigma}
\\
@A{\rm{incl.}}AA @A{\rm{incl.}}AA @ A{\rm{incl.}}AA
\\
\mathscr J @<{\sim}<{\psi_2}< \mathscr K' @>{\sim}>_{\psi_1}>\mathscr K\subset \mathscr U_{\rm{ASD}}
\label{diagram02}
\end{CD}
\end{gather}
(Note again that $\mathscr K = \mathscr U_{\rm{ASD}}$ when $n=3$.) Thus in order to understand the complex structures on
$\mathscr O(-n)$ determined by points on $\mathscr K$, it is enough to understand the complex structures on f\/ibers of
the family over $\mathscr J$.
For this purpose we recall from Sections~\ref{ss:ab} and~\ref{ss:AH} that the Kuranishi family $\mathscr
F_n\to\mathbb{C}^{n-1}$ of $\mathbb{F}_n$ is obtained from the family $\mathscr A_n\to\mathbb{C}^{n-1}$ of af\/f\/ine
bundles by taking a~simultaneous compactif\/ication.
Let $\mathscr L_n:=\mathscr F_n\backslash \mathscr A_n$ be the family of sections at inf\/inity.
We now apply the construction in Section~\ref{ss:sncs} to all f\/ibers of $(\mathscr F_n,\mathscr
L_n)\to\mathbb{C}^{n-1}$ simultaneously.
For this, we need to give an involution $\tau$ on each section to make the variety.
For this purpose we note that since all f\/ibers of $\mathscr F_n\to\mathbb{C}^{n-1}$ have a~common projection to
$\mathbb{CP}^1$ (equipped with the coordinates $u$ and $v$ as before), all f\/ibers of $\mathscr L_n\to\mathbb{C}^{n-1}$
are naturally identif\/ied each other.
Let $\tau_0:\mathbb{CP}^1\to\mathbb{CP}^1$ be an anti-holomorphic involution def\/ined by $\tau_0(u)=-1/\overline u$,
and through the identif\/ication we regard $\tau_0$ as an anti-holomorphic involution which is def\/ined on each f\/iber
of $\mathscr L_n\to\mathbb{C}^{n-1}$.
Then taking this $\tau_0$ as the involution $\tau$ in the construction of Section~\ref{ss:sncs} for any
$(t_1,\dots,t_{n-1}) \in\mathbb{C}^{n-1}$, we obtain a~family of smooth normal crossing surfaces, whose parameter space
is $\mathbb{C}^{n-1}$.
We write this family as
\begin{gather}
\label{rsl1}
\mathscr F_n\cup\overline{\mathscr F}_n\to\mathbb{C}^{n-1}.
\end{gather}
From the construction in Section~\ref{ss:sncs} each f\/iber of this family has a~canonical real structure that
interchanges the two components.
In the notation of Def\/inition~\ref{def:scn}, the f\/iber over the origin of this family is isomorphic to
$\mathbb{F}_n\cup_0\overline{\mathbb{F}}_n$ as a~complex variety with real structure, while on the $l$-th coordinate
axis $\mathbb{C}(t_l)$, f\/ibers are isomorphic to $\mathbb{F}_{n-2l}\cup_l\overline{\mathbb{F}}_{n-2l}$ except over the
origin.

The family~\eqref{rsl1} is in ef\/fect isomorphic to the (abstract) family over $\mathscr J$:

\begin{Lemma}
\label{lemma:isom4}
In a~sufficiently small neighborhood of the origin, the family~\eqref{rsl1} is isomorphic to the family of smooth
normal crossing surfaces over $\mathscr J$.
\end{Lemma}

\begin{proof}
By versality of the Kuranishi family for locally trivial deformations of $D\cup\overline D\simeq\mathbb{F}_n\cup_0
\overline{\mathbb{F}}_n$, we have an induced map, for which we denote by $\alpha$, from a~neighborhood of the origin of
the parameter space $\mathbb{C}^{n-1}$ of~\eqref{rsl1} to that of the last Kuranishi family, such that the pull-back by~$\alpha$ is isomorphic to the family~\eqref{rsl1}.
Though $\alpha$ is not uniquely determined, from naturality, the derivative~$\alpha'(0)$ is nothing but the
Kodaira--Spencer map for~\eqref{rsl1} at the origin.
On the other hand as in~\eqref{isom001} and~\eqref{isom002} we have natural isomorphisms
\begin{gather*}
H^1(\mathbb{F}_n \cup_0 \overline{\mathbb{F}}_n,\Theta) \simeq H^1(\mathbb{F}_n, \Theta_{\mathbb{F}_n,L}) \oplus H^1
(\overline{\mathbb{F}}_n,\Theta_{\overline{\mathbb{F}}_n,\overline L}) \simeq H^1(\mathbb{F}_n, \Theta_{\mathbb{F}_n})
\oplus H^1 (\overline{\mathbb{F}}_n,\Theta_{\overline{\mathbb{F}}_n}),
\end{gather*}
and the composition of the Kodaira--Spencer map $\alpha'(0)$ with these two isomorphisms is an injection onto the real
locus of the last direct sum, because from the construction of the family~\eqref{rsl1}, if we further take the
composition with the projection to the f\/irst factor $H^1(\mathbb{F}_n, \Theta_{\mathbb{F}_n})$ of the last direct sum,
we obviously obtain the Kodaira--Spencer map~\eqref{ks001}, which is an isomorphism.
This means that, in the neighborhood of the origin, the family~\eqref{rsl1} is isomorphic to the real locus of the
Kuranishi family of $\mathbb{F}_n\cup_0\overline{\mathbb{F}}_n$.
Since the isomorphism $D\cup\overline D \simeq \mathbb{F}_n\cup\overline{\mathbb{F}}_n$ respects the real structure, the
last real locus is exactly $\mathscr J$, as desired.
\end{proof}

Now we are able to prove our main result, concerning extendability of the LeBrun metric on $\mathscr O(-n)$ to all
nearby f\/ibers of the above family $\mathscr A_n\to \mathbb{C}^{n-1}$ as an ALE SFK metric:

\begin{Theorem}
\label{thm:main01}
The LeBrun metric on $\mathscr O(-n)$ extends smoothly to all nearby fibers of the family $\mathscr
A_n\to\mathbb{C}^{n-1}$ in~\eqref{af6} of affine $\mathbb{C}$-bundles, as an ALE SFK metric.
\end{Theorem}

\begin{proof}
As in Lemma~\ref{lemma:isom4}, via the induced map $\alpha$, the family $\mathscr F_n\cup\overline{\mathscr
F}_n\to\mathbb{C}^{n-1}$ is isomorphic to the family over $\mathscr J$.
Moreover, as we have already seen, the family over $\mathscr J$ is isomorphic to the $(D\cup\overline D)$-portion of the
deformation of the pair $(Z,D\cup\overline D)$ parameterized by $\mathscr K'$ via the induced map $\psi_2$.
Furthermore, the $Z$-portion of the family of pairs over $\mathscr K'$ is identif\/ied with the family of twistor spaces
over the K\"ahler locus $\mathscr K$ via the map $\psi_1$.
(See the diagram~\eqref{diagram02}.) By the theorem of Pontecorvo~\cite{Pont92}, for any point of $\mathscr K$, the
corresponding twistor space determines an SFK metric on the 4-manifold $\mathscr O(-n)$ up to overall constants.
These SFK metrics can be made to be ALE by multiplying overall constant for each metrics, because the af\/f\/ine bundles
we are considering have a~compactif\/ication by a~$(+n)$-curve.
Via the isomorphisms $\psi_1$, $\psi_2$ and $\alpha$, we conclude that all f\/ibers of the family $\mathscr
A_n\to\mathbb{C}^{n-1}$ admit ALE SFK metrics at least as long as the f\/ibers are suf\/f\/iciently close to the central f\/iber.
The smoothness for the variation of the metrics immediately follows from smoothness for $\psi_1$, $\psi_2$ and $\alpha$.
\end{proof}

We note that as in the above proof, the ALE SFK metrics on all nearby f\/ibers of the central f\/iber are uniquely
determined up to overall constants once we f\/ix the maps $\psi_1$, $\psi_2$ and $\alpha$.

From Theorem~\ref{thm:main01} it is immediate to prove the existence of an ALE SFK metric on any af\/f\/ine
$\mathbb{C}$-bundle over $\mathbb{CP}^1$ of negative degree (see Def\/inition~\ref{def:degree}).
For this, we recall that as we have explained in Section~\ref{ss:ab}, any af\/f\/ine $\mathbb{C}$-bundle over
$\mathbb{CP}^1$ of degree~$-n$ $({\le}-1)$ is a~member of the family $\mathscr A_{n}\to\mathbb{C}^{n-1}$.
Also, from the $\mathbb{C}^*$-action on the total space of $\mathscr A_n$ which is a~lift of the scalar multiplication
on $\mathbb{C}^{n-1}$, any f\/iber over the same line through the origin is mutually biholomorphic except over the
origin.
Thus for any sequence $\{U_{-n}\,|\, n\ge 1\}$ of neighborhoods of the origin in $\mathbb{C}^{n-1}$, the union
$\cup_{n\ge 1} U_{-n}$ contains arbitrary af\/f\/ine $\mathbb{C}$-bundles over $\mathbb{CP}^1$ of negative degree.
Hence by Theorem~\ref{thm:main01} we obtain

\begin{Corollary}\label{cor:exaf}
Any affine $\mathbb{C}$-bundle over $\mathbb{CP}^1$ of negative degree $($see Definition~{\rm \ref{def:degree})} admits an
ALE SFK metric.
\end{Corollary}

Also now it is easy to show the following rigidity result for the LeBrun metric on $\mathscr O(-n)$ when the complex
structure is f\/ixed:

\begin{Proposition}
\label{prop:rigid}
Let $k>0$ and $\Delta$ be a~unit disk in $\mathbb{R}^k$ around the origin, and let $\{g_t\,|\, t\in\Delta\}$ be
a~smooth family of ALE SFK metrics on the complex surface $\mathscr O(-n)$ equipped with the natural complex structure
as a~line bundle.
Assume that $g_0$ is isometric to the LeBrun metric.
Then there exists a~neighborhood $\Delta' \subset \Delta$ of the origin, such that $g_t$ is isometric to the LeBrun
metric up to overall constants for any $t\in\Delta'$.
\end{Proposition}

\begin{proof}
For each $t\in \Delta$ we take a~conformal compactif\/ication $\hat g_t$ of $g_t$ to $\widehat{\mathscr O(-n)}$.
Let $Z_t$ be the twistor space of $\hat g_t$ and $F_t$ be the line bundle $K_{Z_t}^{-1/2}$.
Then by the assumption for complex structure on $\mathscr O(-n)$, for any $t\in \Delta$, the twistor space $Z_t$ has
a~Cartier divisor $D_t\cup \overline D_t(\simeq \mathbb{F}_n\cup_0\overline{\mathbb{F}}_n)$ in the system $|F_t|$
which is biholomorphic to the divisor $D\cup \overline D$ in Proposition~\ref{prop:isom1}.
Hence the family $\{(Z_t,D_t\cup\overline D_t)\,|\, t\in\Delta\}$ gives a~locally trivial deformation of the pair
$(Z_0,D_0 \cup \overline D_0)$ for which the complex structure of $D_0\cup \overline D_0$ does not vary.
By versailty of the Kuranishi family for locally trivial deformations of the pair $(Z_0,D_0\cup \overline D_0)$, there
exist a~neighborhood $\Delta'\subset\Delta$ of the origin and a~smooth map $\varphi:\Delta'\to H^1(\Theta_{Z_0,D_0\cup
\overline D_0})$ which satisf\/ies $\varphi(0) = 0$ and whose pullback of the Kuranishi family is isomorphic to the
original family $\{(Z_t,D_t\cup \overline D_t)\,|\, t\in\Delta'\}$.
But because of the constancy $D_t\cup\overline D_t\simeq D_0\cup \overline D_0$ and the natural isomorphism $H^1 (
\Theta_{Z,  D \cup \overline D} ) \simeq H^1 (\Theta_{ D \cup \overline D} )$ in~\eqref{basicisom1}, $\varphi$ has
to satisfy $\varphi(t)=0$ for any $t\in\Delta'$.
This means that the family $\{Z_t\,|\, t\in \Delta'\}$ itself is a~trivial family.
Hence the conformal classes $[\hat g_t]$ do not vary.
This means the required rigidity of the LeBrun's K\"ahler metric.
\end{proof}

Next we take group actions into account for the moduli problem.
Since LeBrun's metric on~$\mathscr O(-n)$ is ${\rm{U}}(2)$-invariant, its twistor space $Z$ admits
a~${\rm{U}}(2)$-action and the divisor $D\cup\overline D$ is ${\rm{U}}(2)$-invariant.
Hence the cohomology group $H^1(\Theta_Z)$ and $H^1(\Theta_{Z,D\cup\overline D})$ have natural ${\rm{U}}(2)$-actions.
The action on $H^1(\Theta_Z)$ was computed in~\cite{HonCMP2}, and if~$H^1(\Theta_Z)^{\sigma}$ denotes the relevant real
locus, we have, as a~real ${\rm{U}}(2)$-module,
\begin{gather*}
H^1(\Theta_Z)^{\sigma}\simeq S^{n-2}_1\mathbb{C}^2 \oplus S^{n-4}_2\mathbb{C}^2.
\end{gather*}
Here, $S^m_k\mathbb{C}^2:= S^m\mathbb{C}^2\otimes_{\mathbb{C}}\mathbb{C}_k$, where $S^m\mathbb{C}^2$ denotes the $m$-th
symmetric product of the natural representation on $\mathbb{C}^2$, and $\mathbb{C}_k$ is the 1-dimensional
representation obtained by multiplying $(\det)^k$.
(If $m < 0$, $S^m\mathbb{C}^2$ means $0$, and $S^0\mathbb{C}^2$ means the trivial representation on $\mathbb{C}$.) For
the ${\rm{U}}(2)$-action on $H^1(\Theta_{Z,D+\overline D})$, it is immediate from Proposition~\ref{prop:isom1} to derive
the following

\begin{Proposition}
As a~real ${\rm{U}}(2)$-module, we have $H^1(\Theta_{Z,D\cup\overline D})^{\sigma} \simeq S^{n-2}_1\mathbb{C}^2$.
\end{Proposition}

\begin{proof}
By the injectivity of the natural map~\eqref{forget}, $H^1(\Theta_{Z,D\cup\overline D})^{\sigma}$ is naturally
a~subspace of $H^1(\Theta_{Z})^{\sigma}$ which is of course ${\rm{U}}(2)$-invariant.
Since both $S^{n-2}_1\mathbb{C}^2 $ and $S^{n-4}_2\mathbb{C}^2$ are irreducible ${\rm{U}}(2)$-modules,
$H^1(\Theta_{Z,D\cup\overline D})^{\sigma}$ has to coincide with one of these two spaces or the whole space
$S^{n-2}_1\mathbb{C}^2 \oplus S^{n-4}_2\mathbb{C}^2$.
But it has to be $S^{n-2}_1\mathbb{C}^2$ as $\dim H^1(\Theta_{Z,D\cup\overline D})^{\sigma}=2(n-1)$ from
Proposition~\ref{prop:isom1} while $\dim_{\mathbb{R}}S_2^{n-4}\mathbb{C}^2 = 2(n-3)\neq 2(n-1)$.
\end{proof}

Thus connecting the series of the natural isomorphisms
\begin{gather}
\label{isom49}
H^1(\Theta_{Z,D\cup\overline D})^{\sigma}
\; \stackrel{\eqref{basicisom1}}{\simeq}
\;
H^1(\Theta_{D\cup\overline D})^{\sigma}
\; \stackrel{\eqref{isom001}\&\eqref{isom002}}{\simeq}
\;
H^1(\Theta_D)
\;\stackrel{\eqref{ks001}}{\simeq}
\;
\mathbb{C}^{n-1},
\end{gather}
we obtain that the parameter space $\mathbb{C}^{n-1}$ of the family $\mathscr A_n\to\mathbb{C}^{n-1}$ may be
identif\/ied with $S^{n-2}\mathbb{C}^2$ as a~real ${\rm{U}}(2)$-module.

We also have the following result concerning ${\rm{U}}(2)$-action on another cohomology group.
Recall from Proposition~\ref{prop:isom1} that we have $H^1(\Theta_Z(-D-\overline D))\simeq\mathbb{C}$.

\begin{Proposition}
For the LeBrun metric, the natural ${\rm{U}}(2)$-action on $H^1(\Theta_Z(-D-\overline D))\simeq\mathbb{C}$ is trivial.
\end{Proposition}

\begin{proof}
As in~\eqref{LB:key} we have the exact sequence
\begin{gather}
\label{les:10}
0
\longrightarrow
H^0(\Theta_{Z, D\cup\overline D})  \big({\simeq}\, \mathbb{C}^4\big)
\;\stackrel{\iota}{\longrightarrow}
\;
H^0(\Theta_{D\cup\overline D})
\big({\simeq}\, \mathbb{C}^5\big)
\longrightarrow
H^1( \Theta_Z(-D-\overline D))
\longrightarrow
0,
\end{gather}
and hence the space $H^1(\Theta_Z(-D-\overline D))$ $(\simeq\mathbb{C})$ can be identif\/ied with the cokernel of the
injection~$\iota$.
In the space $H^0(\Theta_{D\cup \overline D})$ we have the 2-dimensional subspace generated by the scalar multiplication
on each of $D$ and $\overline D$, where we are viewing these as the compactif\/ication of the line bundle $\mathscr
O(-n)$ for the scalar multiplication.
As the scalar multiplications commute with any element of ${\rm{U}}(2)$, the group ${\rm{U}}(2)$ acts trivially on this
2-dimensional subspace.
Moreover the image of $\iota$ cannot contain this subspace since any real element of $H^0(\Theta_{Z,D\cup\overline D})$
is a~lift of a~conformal Killing f\/ield on $\widehat{\mathscr O(-n)}$ and hence on $D$ and $\overline D$ the vector
f\/ield cannot move independently each other.
This means that the 2-dimensional subspace of $H^0(\Theta_{D\cup\overline D})$ is mapped surjectively to the
1-dimensional space $H^1(\Theta_Z(-D-\overline D))$.
Since the sequence~\eqref{les:10} is ${\rm{U}}(2)$-equivariant, the assertion follows.
\end{proof}

Next we investigate the restrictions of the family of ALE SFK metrics in Theorem~\ref{thm:main01} to the coordinate axes
of $\mathbb{C}^{n-1}$, which provide ${\rm{U}}(1)$-equivariant deformations of the LeBrun metric.
Suppose $n\ge 3$ and for each $1\le l\le n-1$ let $\mathbb{C}(t_l)$ be the $l$-th coordinate axis of $\mathbb{C}^{n-1}$
as before, and let
\begin{gather}\label{1para}
\mathscr A_{n,l}\to \mathbb{C}(t_l)
\qquad
\text{and}
\qquad
\mathscr F_{n,l}\to \mathbb{C}(t_l)
\end{gather}
be the restrictions of the families $\mathscr A_n\to\mathbb{C}^{n-1}$ and $\mathscr F_n\to\mathbb{C}^{n-1}$ respectively
to the $t_l$-axis.
From the $\mathbb{C}^*$-action which is a~lift of the scalar multiplication on $\mathbb{C}^{n-1}$, all f\/ibers of
$\mathscr A_{n,l}\to \mathbb{C}(t_l)$ and $\mathscr F_{n,l}\to \mathbb{C}(t_l)$ are mutually isomorphic except the
central f\/iber for each, and as in~\eqref{zeta'}, f\/ibers of $\mathscr F_{n,l}\to \mathbb{C}(t_l)$ are isomorphic to
$\mathbb{F}_{n-2l}$ except the f\/iber over the origin.
Moreover recalling that~$\mathscr A_n$ is def\/ined by the equation
\begin{gather*}
\zeta_0 = \frac 1{u^n} \zeta_1 + \sum\limits_{l=1}^{n-1} \frac {t_l} {u^l}
\qquad
\text{on}
\quad
U_{01},
\end{gather*}
as in~\eqref{af5}, we obtain that, as an enlargement of the above $\mathbb{C}^*$-action on $\mathscr A_n$ (and $\mathscr
F_n$), the total space of the family $\mathscr A_n$ (and $\mathscr F_n$) carries
a~$(\mathbb{C}^*\times\mathbb{C}^*)$-action def\/ined by
\begin{gather}
\label{taction}
(u,\zeta_0,t_l)
\;\stackrel{(s_1,s_2)}{\longmapsto}
\;
(s_1u,s_2\zeta_0,s_1^ls_2t_l),
\qquad
(s_1,s_2)\in\mathbb{C}^*\times\mathbb{C}^*.
\end{gather}
On the central f\/iber this gives a~$(\mathbb{C}^*\times\mathbb{C}^*$)-action of the toric structure on $\mathscr O(-n)$
or $\mathbb{F}_n$.
Putting $s_1=0$ in~\eqref{taction} gives the original $\mathbb{C}^*$-action on $\mathscr A_n$ and $\mathscr F_n$.
If we def\/ine a~$\mathbb{C}^*$-subgroup $G_l$ of $\mathbb{C}^*\times\mathbb{C}^*$ by
\begin{gather}
\label{G_l}
G_l:=\big\{(s_1,s_2)\in\mathbb{C}^*\times\mathbb{C}^*\,|\, s_1^ls_2=1\big\},
\end{gather}
which acts trivially on $\mathbb{C}(t_l)$, then the two families~\eqref{1para} may be regarded as $G_l$-equivariant
deformations of $\mathscr O(-n)$ and $\mathbb{F}_n$ respectively.
Then basically by restricting the family of ALE SFK metrics in Theorem~\ref{thm:main01} to the coordinate axis
$\mathbb{C}(t_l)$, we obtain the following result about existence of ${\rm{U}}(1)$-equivariant deformations of the
LeBrun metric:

\begin{Proposition}
\label{prop:invm}
Let $n\ge 3$ and $l\ge 1$ be integers satisfying $n-2l\ge 0$, and $\mathscr A_{n,l} \to\mathbb{C}(t_l)$ be the
$1$-parameter deformation of $\mathscr O(-n)$ to an affine $\mathbb{C}$-bundle as above.
Then the LeBrun metric on the central fiber $\mathscr O(-n)$ extends smoothly to any nearby fibers of $\mathscr
A_{n,l}\to\mathbb{C}(t_l)$ preserving not only the ALE SFK property but also a~${\rm{U}}(1)$-action, where ${\rm{U}}(1)$
is the compact torus of the stabilizer subgroup $G_l(\simeq\mathbb{C}^*)$ of the axis $\mathbb{C}(t_l)$ defined in~\eqref{G_l}.
\end{Proposition}

\begin{proof}
From the proof of Theorem~\ref{thm:main01}, we just need to show that the isomorphism between the K\"ahler locus
$\mathscr K$ and the parameter space of the family $\mathscr F_n\cup \overline{\mathscr F}_n\to\mathbb{C}^{n-1}$ around
the origin can be taken to be $T^2$-equivariant, where $T^2$ is the standard maximal torus of ${\rm{U}}(2)$ consisting
of diagonal matrices.
For this it is enough to see that the induced maps $\psi_1$, $\psi_2$ and $\alpha$, which were used to identify the
relevant families in Section~\ref{ss:DLm}, can be taken to be $T^2$-equivariant.
This holds for $\psi_1$ and $\psi_2$ since $Z$, $D\cup\overline D$ and the pair $(Z,D\cup\overline D)$ are
${\rm{U}}(2)$-invariant.
For $T^2$-equivariance of the remaining map $\alpha$, it is enough to see that the total space of the family $\mathscr
F_n\cup\overline{\mathscr F}_n\to\mathbb{C}^{n-1}$ has a~$T^2$-action whose restriction to the central f\/iber is
identif\/ied with the $T^2$-action on $D\cup\overline D$.
Each component of $\mathscr F_n\cup\overline{\mathscr F}_n$ has a~$T^2$-action which is the restriction of the
$(\mathbb{C}^*\times\mathbb{C}^*)$-action given in~\eqref{taction} to the maximal torus, and these are clearly
identif\/ied with the $T^2$-actions on $D$ and $\overline D$ respectively.
So to complete the proof we just need to see that the gluing map which was used for making the family $\mathscr
F_n\cup\overline{\mathscr F}_n$ is $T^2$-equivariant.
But this is immediate if we notice that the anti-podal map $\tau_0$ commutes with the $T^2$-actions.
\end{proof}

\begin{Remark}
The parameter space $\mathbb{C}(t_l)$ of the above family of the metric is of course real 2-dimensional, but the family
is in ef\/fect real 1-dimensional by the following reason.
As in the above proof, the axis $\mathbb{C}(t_l)$ is naturally identif\/ied (via the isomorphisms in~\eqref{isom49})
with a~$T^2$-invariant subspace of $H^1(\Theta_{Z,D\cup \overline D})^{\sigma}$ on which the subgroup $G_l$ acts
trivially.
This $T^2$-action has clearly real 1-dimensional orbits, and along each orbit the complex structure of the pairs of
twistor spaces and the divisors are constant.
We also note that although the group $T_{\mathbb{C}}$ (and $T_{\mathbb{C}}/G_l$) is acting on the axis
$\mathbb{C}(t_l)\subset\mathbb{C}^{n-1}$, the corresponding subspace is {\em not} $T_{\mathbb{C}}$- (nor
$T_{\mathbb{C}}/G_l$-) invariant, because the isomorphism in~\eqref{isom49} are not $T_{\mathbb{C}}$-equivariant and
just $T^2$-equivariant.
\end{Remark}

By Corollary~\ref{cor:unique} and Proposition~\ref{prop:invm}, we obtain the following

\begin{Proposition}
Let $n\ge 0$, $l\ge 0$, and let $L\in |\Gamma_0+(n+l)f|$ be any $\mathbb{C}^*$-invariant section of
$\mathbb{F}_n\to\mathbb{CP}^1$, where $\mathbb{C}^*$ acts non-trivially on~$L$.
Then the complement $\mathbb{F}_n\backslash L$ admits a~${\rm{U}}(1)$-invariant ALE SFK metric, where ${\rm{U}}(1)$ is
the compact torus of $\mathbb{C}^*$.
\end{Proposition}

\begin{proof}
By Corollary~\ref{cor:unique}, the complex structure of the pair $(\mathbb{F}_n,L)$ satisfying the properties in the
proposition is uniquely determined from $n$ and $l$.
If $l=0$, we have $\mathbb{F}_n\backslash L\simeq\mathscr O(-n)$ and the existence of the metric on it is guaranteed by
the original LeBrun metric on $\mathscr O(-n)$.
If $l>0$, $\mathbb{F}_n\backslash L$ is biholomorphic to general f\/ibers of the 1-parameter family $\mathscr
A_{n+2l,l}\to\mathbb{C}(t_l)$, and the existence of the metric is guaranteed by Proposition~\ref{prop:invm}, as long as
$n+2l\ge 3$.
The situation where $n+2l\ge 3$ does not hold is only the case $(n,l) = (0,1)$.
But the existence of an ALE SFK metric on $\mathbb{F}_0\backslash L$ ($L\in|\mathscr O(1,1)|$) is guaranteed by the
Eguchi--Hanson metric.
\end{proof}

When $l>1$, if~$L$ is the $\mathbb{C}^*$-invariant section as in the above proposition, then we have
${\rm{Aut}}_0(\mathbb{F}_n,L) \simeq\mathbb{C}^*$ by Proposition~\ref{prop:l2}.
But when $l=1$, we have ${\rm{Aut}}_0(\mathbb{F}_{n}, L)\simeq \rm{Af}(\mathbb{C})$ by Proposition~\ref{prop:nonred}.
Thus the af\/f\/ine surface $\mathbb{F}_n\backslash L$ admits an ALE SFK metric with an ef\/fective ${\rm{U}}(1)$-action
but ${\rm{Aut}}_0(\mathbb{F}_{n}, L)$ is not reductive.
But we do not know whether the holomorphic transformation group of the complex surface $\mathbb{F}_n\backslash L$ itself
is reductive, nor even whether it is of f\/inite dimensional.
In this respect, for the surface $\mathscr O(-n)$, the holomorphic automorphism group is known to be {\em not} of
f\/inite dimensional~\cite[Remark 2.20]{KPZ}.

\subsection{Deformations of the metrics on the af\/f\/ine bundles}

In the last subsection we have obtained ALE SFK metrics on af\/f\/ine $\mathbb{C}$-bundles over $\mathbb{CP}^1$ as small
deformations of the LeBrun metric on $\mathscr O(-n)$.
In this subsection we in turn investigate small deformations of these metrics on af\/f\/ine bundles, which again
preserve ALE SFK property.
So let $\mathscr A_n\to\mathbb{C}^{n-1}$ ($n\ge 3$) be the family~\eqref{af6} of af\/f\/ine bundles over $\mathbb{CP}^1$
as before, and for each $t\in\mathbb{C}^{n-1}$ write $A_t$ for the af\/f\/ine bundle lying over $t$.
By Theorem~\ref{thm:main01}, if $t$ is suf\/f\/iciently close to the origin, the af\/f\/ine bundle $A_t$ admits an ALE
SFK metric.
We write $g_t$ for this metric.
We recall that these metrics are uniquely determined up to overall constants by the (natural but non-unique) maps
$\psi_1$, $\psi_2$ and $\alpha$.
Then these metrics satisfy the following property:

\begin{Proposition}\label{prop:eg}
If $t\neq 0$, the metric $g_t$ on the affine bundle $A_t$ is a~member of a~non-trivial $($see below$)$, real $1$-parameter
family of ALE SFK metrics, in which the complex structure on $A_t$ does not deform.
\end{Proposition}

Here, `non-trivial family' means that the complex structures of the corresponding 1-parameter family of twistor spaces
actually vary as the parameter moves.
Thus the situation is in contrast with the LeBrun metric on the line bundle $\mathscr O(-n)$, for which the metric
cannot be deformed as an ALE SFK metric when the complex structure is f\/ixed (see Proposition~\ref{prop:rigid}).

\begin{proof}[Proof of Proposition~\ref{prop:eg}] Let $Z_t$ be the twistor space of a~conformal compactif\/ication of~$g_t$, and
let~$D_t$ be the divisor determined by the complex structure of~$A_t$.
The sum $D_t+\overline D_t$ is a~Cartier divisor belonging to $|K_{Z_t}^{-1/2}|$.
We f\/irst show that the natural map
\begin{gather}
\label{surj}
H^0(\Theta_{Z_t,D_t\cup\overline D_t})
\longrightarrow
H^0(\Theta_{D_t\cup\overline D_t})
\end{gather}
is surjective as long as $t\neq 0$ and $t$ is suf\/f\/iciently close to $0$.
(Note that if $t=0$ this is not surjective as in the sequence~\eqref{les:1}.) This is trivially satisf\/ied if
$H^0(\Theta_{D_t\cup\overline D_t}) = 0$.
If $H^0(\Theta_{D_t\cup\overline D_t}) \neq 0$, we clearly have $H^0(\Theta_{D_t,L_t}) \neq 0$, where
$L_t=D_t\cap\overline D_t$.
From Propositions~\ref{prop:l1},~\ref{prop:nonred} and~\ref{prop:l2}, this is the case only when $(D_t,L_t)$ is
a~$\mathbb{C}^*$-invariant pair.
Again by the same propositions, the complex structure of such a~pair is unique once the two integers $m\ge 0$ and $l\ge
1$ are f\/ixed for which $D_t\simeq\mathbb{F}_m$ and $L\in |\Gamma_0 + (m+l)f|$ hold.
This means that we have $D_t\cup \overline D_t\simeq \mathbb{F}_{n-2l}\cup_l\overline{\mathbb{F}}_{n-2l}$ for some $l\ge
1$ (as $t\neq 0$) in the notation of Def\/inition~\ref{def:scn}.
For the central f\/iber $(Z_0,D_0\cup\overline D_0)$, as in~\eqref{basicisom1}, we have a~natural isomorphism
\begin{gather}
\label{edisom}
H^1( \Theta_{Z_0,D_0\cup\overline D_0}) \simeq H^1(\Theta_{D_0\cup\overline D_0}),
\end{gather}
and this is ${\rm{U}}(2)$-equivariant.
Since $H^2( \Theta_{Z_0,D_0\cup\overline D_0}) = H^2(\Theta_{D_0\cup\overline D_0})=0$ as
in~\eqref{van1},~\eqref{edisom} means that, for any subgroup $G\subset {\rm{U}}(2)$, $G$-action on the surface
$D_0\cup\overline D_0$ extends to $D_t\cup\overline D_t$ for suf\/f\/iciently small $t$ if and only if the $G$-action on
the pair $(Z_0,D_0\cup\overline D_0)$ extends to $(Z_t,D_t\cup\overline D_t)$.
Moreover the restriction of the family $\mathscr F_n\cup\overline{\mathscr F}_n\to\mathbb{C}^{n-1}$ given
in~\eqref{rsl1} to the coordinate axis~$\mathbb{C}(t_l)$ actually connects $D_0\cup\overline D_0$ and $D_t\cup\overline
D_t$ in a~${\rm{U}}(1)$-equivariant way, where ${\rm{U}}(1)$ is the compact torus of~$G_l$ def\/ined in~\eqref{G_l}.
Thus the ${\rm{U}}(1)$-action on $(Z_0,D_0\cup\overline D_0)$ actually extends to $(Z_t,D_t\cup\overline D_t)$.
This means that the map~\eqref{surj} is surjective for suf\/f\/iciently small $t$.

By the upper semi-continuity and the invariance of the Euler characteristic under deformation, we have
\begin{gather*}
H^i(Z_t,\Theta_{Z_t}(-D_t-\overline D_t)) =0
\qquad
\text{if}
\quad
i\neq 1,
\qquad
\text{and}
\qquad
H^1(Z_t,\Theta_{Z_t}(-D_t-\overline D_t)) \simeq\mathbb{C},
\end{gather*}
because these are true for the case of the LeBrun metric as in~\eqref{cohomdim1}.
Therefore from the exact sequence~\eqref{ses:1} (which remains obviously valid for $(Z_t,D_t\cup\overline D_t)$) and the
surjectivity of the map~\eqref{surj}, we obtain the short exact sequence
\begin{gather}
\label{ses:10}
0
\longrightarrow
H^1(\Theta_{Z_t}(-D_t-\overline D_t))
(\simeq\mathbb{C})
\;\stackrel{\alpha}{\longrightarrow}
\;
H^1(\Theta_{Z_t, D_t\cup\overline D_t})
\;\stackrel{\beta}{\longrightarrow}
\;
H^1(\Theta_{D_t\cup \overline D_t})
\longrightarrow
0.
\end{gather}
Then as we have $H^2(\Theta_{Z_t,D_t+\overline D_t}) = 0$ by~\eqref{van1} and the upper semi-continuity again, the
parameter space of the Kuranishi family for locally trivial deformations of the pair $(Z_t,D_t\cup\overline D_t)$ can be
regarded as a~small disk (for which we denote by $\Delta_1$) about the origin in $H^1(\Theta_{Z_t,D_t\cup\overline
D_t})$.
Similarly as we have $H^2(\Theta_{D_t\cup\overline D_t})=0$ from~\eqref{van1}, the Kuranishi family for $D_t\cup
\overline D_t$ may be regarded as a~small disk (for which we denote by $\Delta_2$) about the origin in
$H^1(\Theta_{D_t\cup\overline D_t})$.
If $\varphi:\Delta'_1\to\Delta_2$ denotes a~holomorphic map from a~possibly smaller disk $\Delta'_1\subset\Delta_1$
which is induced by the versality of the Kuranishi family for $D_t\cup\overline D_t$, then we naturally have
$\varphi'(0) = \beta$.
Hence from the surjectivity of $\beta$ in~\eqref{ses:10}, $\varphi$ is a~submersion at least in a~neighborhood of the
origin, and therefore $\varphi^{-1}(0)$ is non-singular and 1-dimensional near the origin, at which the tangent space is
exactly the line $\Image(\alpha)$.
Thus the Kuranishi family for locally trivial deformations of the pair $(Z_t,D_t\cup \overline D_t)$ contains
a~(complex) 1-parameter subfamily over which the complex structure of $D_t\cup\overline D_t$ is constant.
By restricting to the real locus of $\varphi^{-1}(0)$, we obtain the real 1-parameter family of deformation of the pair
$(Z_t,D_t\cup \overline D_t)$ for which the complex structure of $D_t\cup\overline D_t$ is constant.
By the same reason as in the proof of Theorem~\ref{thm:main01}, the ASD structure determined by $(Z_t,D_t\cup\overline
D_t)$ has a~K\"ahler representative which is ALE at inf\/inity.

It remains to show that the 1-parameter family of ALE SFK metrics on the 4-manifold $\mathscr O(-n)$ thus obtained is
non-trivial in the sense explained right after the proposition.
For this, it suf\/f\/ices to see that the Kodaira--Spencer class of the 1-parameter family of $Z_t$ is non-zero, because
this means that the complex structure of $Z_t$ actually deforms.
From the above argument the Kodaira--Spencer class of the deformation of the pair $(Z_t,D_t\cup\overline D_t)$ is
non-zero (belonging to $\Image(\alpha)$), and the genuine Kodaira--Spencer class (belonging to $H^1(\Theta_{Z_t})$) is
obtained from this by sending it under the natural map $H^1(\Theta_{Z_t,D_t\cup\overline D_t})\to H^1(\Theta_{Z_t})$.

We show that the last map is injective in the same way to the last part in the proof of Proposition~\ref{prop:isom1}.
Exactly by the same argument as to deduce~\eqref{N''}, for the cokernel sheaf $N'_t$ of the natural injection
$\Theta_{Z_t,D_t\cup\overline D_t} \to \Theta_{Z_t}$, we obtain
\begin{gather*}
N_t'|_{D_t}\simeq -K_{D_t}-\mathscr O_{D_t}(2L_t),
\qquad
N_t'|_{\overline D_t}\simeq -K_{\overline D_t}-\mathscr O_{\overline D_t}(2L_t).
\end{gather*}
Now as $t\neq 0$ we can write $D_t\simeq\mathbb{F}_{n-2l}$ for some $l$ satisfying $1\le l\le n-1$, and we have $L_t\in
|\Gamma_0+(n-l)f|$.
Hence by the same computation to deduce~\eqref{N'} we obtain
\begin{gather*}
-K_{D_t}-\mathscr O_{D_t}(2L_t) \simeq 2\Gamma_0+(n-2l+2)f -2\{\Gamma_0 + (n-l)f\}
\simeq -(n-2)f.
\end{gather*}
Thus as $n\!>\!2$ we obtain $H^0(-K_{D_t}\!-\mathscr O_{D_t}(2L_t))\!=\!0$.
With reality, this means \mbox{$H^0(D_t\!\cup\! \overline D_t,N'_t)\!=\!0$}.
Hence the injectivity of the map $H^1(\Theta_{Z_t,D_t\cup\overline D_t})\to H^1(\Theta_{Z_t})$ follows, and we obtained
the assertion of the proposition.
\end{proof}

\begin{proof}[Proof of Proposition~\ref{prop:var}]
Let $\mathscr Z\to B$ be the family of twistor spaces associated to ALE SFK
metrics on nearby f\/ibers for the central f\/iber of the family $\mathscr A_n\to \mathbb{C}^{n-1}$ of af\/f\/ine
$\mathbb{C}$-bundle over~$\mathbb{CP}^1$.
For $t\in B$ we write $Z_t$ for the twistor space over the af\/f\/ine bundle $A_t$.
Since the family $\mathscr Z\to B$ is versal at the origin, it is also versal at $t\in B$ as long as $t$ is
suf\/f\/iciently close to the origin.
Therefore for a~real 1-parameter family of twistor spaces associated to the family of ALE SFK metrics obtained in
Proposition~\ref{prop:eg}, there exists an induced map from the parameter space of the last family to~$B$.
Then if we take the image of the map as the arc, it clearly satisf\/ies the property of the proposition.
\end{proof}

\subsection*{Acknowledgements}

I would like to thank Jef\/f Viaclovsky for numerous useful discussion from which this work has originated.
Also I would like to thank Akira Fujiki for invaluable comments on deformations of line bundles.
This work was supported by JSPS KAKENHI Grant Number 24540061.

\pdfbookmark[1]{References}{ref}
\LastPageEnding

\end{document}